\documentclass[11pt,a4paper]{article}
\usepackage{amsfonts,amsmath,amssymb,accents,bm,amsthm}
\usepackage{verbatim}

\usepackage{hyperref}
\usepackage[normalem]{ulem}

\usepackage{tikz} 
\usetikzlibrary{math,calc,arrows,positioning,fit,petri,external}

\usepackage{graphicx}

\newcommand{\dis}{\displaystyle}
\newcommand{\txt}{\textstyle}

\newcommand{\noi}{\noindent}
\newcommand{\blok}{\hspace{-1pt}\tikz[scale=1,baseline=-4pt]{\draw[line width=5pt] (0,0) -- ++ (0.4,0);}\hspace{1pt}}
\newcommand{\halmos}{\rule{1ex}{1.4ex}}
\newcommand{\QED}{\nopagebreak{\hspace*{\fill}$\halmos$\medskip}}
\newcommand{\med}{\medskip}
\newcommand{\quand}{\quad\mbox{and}\quad}

\newtheoremstyle{mythm}
  {}
  {}
  {\itshape}
  {}
  {\bfseries}
  {}
  {.5em}
  {#1 #2 \thmnote{(#3)}}

\theoremstyle{mythm}
\newtheorem{theorem}{Theorem}
\newtheorem{proposition}[theorem]{Proposition}
\newtheorem{lemma}[theorem]{Lemma}

\newcommand{\bt}{\begin{theorem}}
\newcommand{\et}{\end{theorem}}
\newcommand{\bl}{\begin{lemma}}
\newcommand{\el}{\end{lemma}}
\newcommand{\bp}{\begin{proposition}}
\newcommand{\ep}{\end{proposition}}

\newenvironment{Proof}[1][]{\noi\textbf{Proof #1}}{\QED}
\newcommand{\bpro}{\begin{Proof}}
\newcommand{\epro}{\end{Proof}}

\newcommand{\be}{\begin{equation}}
\newcommand{\ee}{\end{equation}}
\newcommand{\ba}{\begin{array}}
\newcommand{\ea}{\end{array}}
\newcommand{\bac}{\begin{array}{r@{\,}c@{\,}l}}
\newcommand{\bc}{\be\begin{array}{r@{\,}c@{\,}l}}
\newcommand{\ec}{\end{array}\ee}

\newcommand{\ga}{\gamma}
\newcommand{\Ga}{\Gamma}
\newcommand{\de}{\delta}
\newcommand{\De}{\Delta}

\newcommand{\la}{\lambda}
\newcommand{\La}{\Lambda}
\newcommand{\sig}{\sigma}

\newcommand{\om}{\omega}

\newcommand{\Ci}{{\cal C}}
\newcommand{\Di}{{\cal D}}

\newcommand{\Gi}{{\cal G}}

\newcommand{\Li}{{\cal L}}

\newcommand{\Pc}{{\cal P}}

\newcommand{\Ri}{{\cal R}}

\newcommand{\E}{{\mathbb E}}

\newcommand{\N}{{\mathbb N}}

\renewcommand{\P}{{\mathbb P}}

\newcommand{\R}{{\mathbb R}}

\newcommand{\Xb}{{\mathbf X}}
\newcommand{\Yb}{{\mathbf Y}}

\newcommand{\Lab}{{\mathbf\La}}
\newcommand{\mb}{{\mathbf m}}

\newcommand{\ob}{{\mathbf o}}
\newcommand{\rb}{{\mathbf r}}
\newcommand{\idb}{{\mathbf{id}}}
\newcommand{\ibf}{\mathbf{i}}
\newcommand{\jbf}{\mathbf{j}}

\newcommand{\omb}{{\bm\omega}}

\newcommand{\Gk}{\mathfrak{G}}

\newcommand{\desd}{\ensuremath{\Leftrightarrow}}
\newcommand{\volgt}{\ensuremath{\Rightarrow}}
\newcommand{\up}{\uparrow}
\newcommand{\down}{\downarrow}
\newcommand{\sub}{\subset}
\newcommand{\beh}{\backslash}

\newcommand{\ti}{\tilde}

\newcommand{\dgg}{\dagger}
\newcommand{\ov}{\overline}
\newcommand{\un}{\underline}

\newcommand{\cn}{\colon}

\newcommand{\half}{{[0,\infty)}}

\begin{document}

\makeatletter\@addtoreset{equation}{section}
\makeatother\def\theequation{\thesection.\arabic{equation}}

\renewcommand{\labelenumi}{{\rm (\roman{enumi})}}
\renewcommand{\theenumi}{\roman{enumi}}

\title{Percolation representations of additive particle systems}
\author{Jan~M.~Swart\footnote{The Czech Academy of Sciences, Institute of Information Theory and Automation, Pod vod\'arenskou v\v{e}\v{z}\'i~4, 18200 Praha 8, Czech republic. swart@utia.cas.cz}} 


\maketitle

\begin{abstract}\noi
It is well-known that additive interacting particle systems with a local state space of cardinality two have a percolation representation in terms of open paths in a graphical representation. In this paper, it is shown how such a percolation representation can be constructed more generally when the local state space is a finite distributive lattice. The theory is demonstrated on Krone's two-stage contact process.
\end{abstract}
\vspace{.5cm}

\noi

\noi


\section{Introduction and statement of the results}

\subsection{Introduction}\label{S:intro}

An \emph{interacting particle system} is a contintuous-time Markov process $(X_t)_{t\geq 0}$ taking values in the space $S^\La$ of all functions from a countable set $\La$, called the \emph{grid}\footnote{This is often called the \emph{lattice}, but we will need that word in its order theoretic meaning.} into a finite set $S$, called the \emph{local state space}. Its generator $G$ can usually be written in the form
\be\label{Gf}
Gf(x):=\sum_{m\in\Gi}r_m\big\{f\big(m(x)\big)-f\big(x\big)\big\},
\ee
where $\Gi$ is a countable collection of maps $m\cn S^\La\to S^\La$ and $(r_m)_{m\in\Gi}$ are nonnegative rates. The maps $m\in\Gi$ are \emph{local}, in the sense that for each $m$ there is a finite set $\Ga\sub\La$ such that $m$ leaves the values of a configuration $x\in S^\La$ outside $\Ga$ unchanged and there is a map $m'\cn S^\Ga\to S^\Ga$ that determines how $m$ changes $x$ inside $\Ga$.


Formula (\ref{Gf}) has the following interpretation: for each $m\in\Gi$, there is an independent Poisson point process on $\half$  with intensity $r_m$ that determines the times $t$ when the Markov process $(X_t)_{t\geq 0}$ jumps from its previous state $X_{t-}$ to its new state $m(X_t)$. It is often useful to extend these Poisson processes from $\half$ to the entire real line $\R$. In this case, under suitable assumptions on the rates $r_m$ (to be discussed below), it can be shown that for each pair of times $s,t\in\R$ with $s\leq t$, there is a random map $\Xb_{s,t}\cn S^\La\to S^\La$ that determines how the process evolves between the times $s$ and $t$. More precisely, starting with the state $X_s$ of the process at time $s$ and locally changing the configuration by applying maps $m$ from $\Gi$ at the correct times as determined by the Poisson point processes with intensity $r_m$, one obtains the state $X_t=\Xb_{s,t}(X_s)$ of the process at time $t$. In this context, the collection of Poisson point processes is called a \emph{graphical representation}, and the collection of random maps $(\Xb_{s,t})_{s\leq t}$ is called the \emph{stochastic flow} constructed from this graphical representation.

In the special case that $S=\{0,1\}$, a map $m\cn S^\La\to S^\La$ is called \emph{additive} if it maps the all zero configuration $\un 0$ into itself and it maps the pointwise maximum $x\vee y$ of two configurations into $m(x)\vee m(y)$. If the local maps $m\in\Gi$ are all additive, then so are the random maps $\Xb_{s,t}$ $(s\leq t)$. Some of the most well-known interacting particle systems such as the contact process, the voter model, and the symmetric exclusion process are additive. Additive local maps can be represented in terms of arrows and blocking symbols (to be discussed below) and the random maps $\Xb_{s,t}$ then have an interpretation in terms of open paths. This allows one to understand the graphical representation of an additive interacting particle system as a form of oriented percolation. These sort of percolation representations of additive interacting particle systems originate from the work of Harris \cite{Har78} and Griffeath \cite{Gri79} and they have proved to be very useful. Tracing open paths backwards in time, one can use a percolation representation to define a second interacting particle system that is dual to the original one. For example, the additive dual of the voter model is a system of coalescing random walks while the contact process and symmetric exclusion process are self-dual.

Additive interacting particle systems can be defined more generally when the local state space $S$ is a finite lattice. An example of an additive interacting particle system with local state space $S=\{0,1,2\}$ is the \emph{two-stage contact process} introduced by Krone \cite{Kro99}. As some authors have realised \cite{Fox16,SS18}, the duality theory for additive interacting particle systems carries over to the case when $S$ is any finite lattice. In \cite{SS18}, it was moreover claimed that as long as the lattice $S$ is distributive,\footnote{Recall that a lattice $S$ is \emph{distributive} if it satisfies one, and hence both of the equivalent conditions (i) $a\wedge(b\vee c)=(a\wedge b)\vee(a\wedge c)$ $(a,b,c\in S)$, and (ii) $a\vee(b\wedge c)=(a\vee b)\wedge(a\vee c)$ $(a,b,c\in S)$.} it is also possible to give a percolation representation. On closer inspection, the construction in \cite{SS18} works only for finite grids $\La$, which is unsatisfactory in practice. In the present paper, we correct this by giving a different construction that shows that each additive interacting particle system whose local state space is a distributive lattice has a percolation representation in terms of open paths in a graphical representation involving arrows and blocking symbols.

\subsection{Interacting particle systems}\label{S:IPS}

In this subsection, we recall some basic facts about the construction of interacting particle systems from a graphical representation. We follow the approach of \cite{Swa26}. Let $S$ and $\La$ be a finite and countable set as before. We equip the space $S^\La$ of all functions $x\cn\La\to S$ with the product topology. If $T$ is a finite set and $f\cn S^\La\to T$ is a function, then we say that a grid point $i\in\La$ is \emph{$f$-relevant} if changing a configuration $x$ in the grid point $i$ can potentially change the value of $f(x)$, that is
\be\label{relev}
\exists x,y\in S^\La\mbox{ s.t.\ }f(x)\neq f(y)\mbox{ and }x(j)=y(j)\ \forall j\neq i.
\ee
We let $\Ri(f):=\{i\in\La:i\mbox{ is $f$-relevant}\}$ denote the set of grid points that are relevant for $f$. By \cite[Lemma~2.30]{Swa26}, a function $f\cn S^\La\to T$ is continuous if and only if it satisfies the conditions
\begin{enumerate}
\item $\Ri(f)$ is finite,
\item if $x,y\in S^\La$ satisfy $x(i)=y(i)$ for all $i\in\Ri(f)$, then $f(x)=f(y)$.
\end{enumerate}
Similarly, for a map $m\cn S^\La\to S^\La$ and $j\in\La$, we define $m[j]\cn S^\La\to S$ by $m[j](x):=m(x)(j)$, and we set
\be\ba{c}\label{DR}
\dis\Di(m):=\big\{i\in\La:\exists x\in S^\La\mbox{ s.t.\ }m(x)(i)\neq x(i)\big\},\\[5pt]
\dis\Ri(m):=\big\{(i,j)\in\La^2:i\mbox{ is $m[j]$-relevant}\big\},\\[5pt]
\dis\Ri^\down_j(m):=\big\{i\in\La:(i,j)\in\Ri(m)\big\},\quad
\Ri^\up_i(m):=\big\{j\in\La:(i,j)\in\Ri(m)\big\}.
\ec
By definition, a \emph{local} map $m\cn S^\La\to S^\La$ is a continuous map for which $\Di(m)$ is finite. It is not difficult to prove that this is equivalent to the statement that there exists a finite $\Ga\sub\La$ and a map $m'\cn S^\Ga\to S^\Ga$ such that $m(x)(i)=m'(x|_\Ga)(i)$ for $i\in\Ga$ and $m(x)(i)=x(i)$ for $i\in\La\beh\Ga$, where $x|_\Ga$ denotes the restriction of $x$ to $\Ga$; see \cite[Exercise~2.33]{Swa26}. Let $\Gi$ be a countable collection of local maps and let $(r_m)_{m\in\Gi}$ be nonnegative rates. We will need the conditions
\be\ba{c}\label{sum}
\dis{\rm(i)}\ \sup_{i\in\La}\sum_{m\in\Gi}r_m1_{\Di(m)}(i)<\infty,\\[5pt]
\dis{\rm(ii)}\ \sup_{i\in\La}\sum_{m\in\Gi}r_m\big|\Ri^\down_i(m)\beh\{i\}\big|<\infty,\quad
{\rm(iii)}\ \sup_{i\in\La}\sum_{m\in\Gi}r_m\big|\Ri^\up_i(m)\beh\{i\}\big|<\infty.
\ec
By \cite[Thms 4.14 and Thm~4.29]{Swa26}, conditions (\ref{sum}) (i) and (ii) guarantee that the Markov process with generator $G$ as in (\ref{Gf}) is well-defined and can be constructed from a graphical representation. If moreover (\ref{sum}) (iii) holds, then by \cite[Prop~4.24]{Swa26}, finite perturbations of the initial state lead to finite perturbations at later times. We now describe these results more precisely.

Let $\om$ be a Poisson point set on $\Gi\times\R$ with intensity $\rho$ given by
\be\label{rho}
\rho\big(\{m\}\times[s,t]\big):=r_m(t-s)\qquad(m\in\Gi,\ s\leq t).
\ee
Recall that a function $f$, defined on a real interval $I$ and taking values in a metrisable space, is called \emph{cadlag} (from the French continue \`a droit, limite \`a gauche) if it is right-continuous with left limits $f(t-):=\lim_{s\up t}f(s)$ in each point $t\in I$ that can be approximated from the left. By \cite[Thm~4.14]{Swa26}, under the conditions (\ref{sum}) (i) and (ii), almost surely, for each $x\in S^\La$ and $s\in\R$, there exists a unique cadlag function $X^{s,x}\cn[s,\infty)\to S^\La$ with initial state $X^{s,x}_s=x$ that satisfies the evolution equation
\be\label{evol}
X^{s,x}_t:=\left\{\ba{ll}
\dis m(X^{s,x}_{t-})\quad&\mbox{if }(m,t)\in\om,\\[5pt]
X^{s,x}_{t-}\quad&\mbox{otherwise,}\ea\right.\qquad(s<t).
\ee
Moreover, setting $\Xb_{s,t}(x):=X^{s,x}_t$ $(s\leq t,\ x\in S^\La)$ defines a collection $(\Xb_{s,t})_{s\leq t}$ of continuous maps from $S^\La$ into itself such that
\be\label{oldflow}
\Xb_{s,s}=1\quand\Xb_{t,u}\circ\Xb_{s,t}=\Xb_{s,u}\qquad(s\leq t\leq u),
\ee
where $1$ denotes the identity map. We call $(\Xb_{s,t})_{s\leq t}$ the \emph{stochastic flow} defined by the \emph{graphical representation} $\om$. Since restrictions of a Poisson point process to disjoint sets are independent, the stochastic flow $(\Xb_{s,t})_{s\leq t}$ has \emph{independent increments} in the sense that for each $t_0<\cdots<t_n$, the random variables $(\Xb_{t_{k-1},t_k})_{1\leq k\leq n}$ are independent. It is also \emph{stationary} in the sense that $(\Xb_{s,t})_{s\leq t}$ is equally distributed with $(\Xb_{r+s,r+t})_{s\leq t}$ for each $r\in\R$. By \cite[Thm~4.14]{Swa26}, setting
\be
P_t(x,\,\cdot\,):=\P\big[\Xb_{0,t}(x)\in\,\cdot\,\big]\qquad(t\geq 0,\ x\in S^\La)
\ee
defines a Feller semigroup $(P_t)_{t\geq 0}$ such that if $s\in\R$ and $X_0$ is an $S^\La$-valued random variable with law $\mu$, independent of $\om$, then setting
\be\label{XX}
X_t:=\Xb_{s,s+t}(X_0)\qquad(t\geq 0)
\ee
defines a Feller process $(X_t)_{t\geq 0}$ with semigroup $(P_t)_{t\geq 0}$. By \cite[Thm~4.29]{Swa26}, the generator of $(P_t)_{t\geq 0}$ is the closure of the operator $G$ from (\ref{Gf}), which is initially defined only for functions depending on finitely many coordinates.

By \cite[Prop~4.24]{Swa26}, if in addition to (\ref{sum}) (i) and (ii) also condition (\ref{sum}) (iii) is satisfied, then almost surely, for each $s\leq t$ and $x,y\in S^\La$,
\be\label{finpert}
\big|\big\{i\in\La:x(i)\neq y(i)\big\}\big|<\infty\quad\mbox{implies}\quad
\big|\big\{i\in\La:\Xb_{s,t}(x)(i)\neq\Xb_{s,t}(y)(i)\big\}\big|<\infty,
\ee
which says that under the stochastic flow, finite perturbations of the initial state lead to finite perturbations at later times. 

For later use, we note that each graphical representation $\om$ can also be used to define a \emph{backward stochastic flow}. Indeed, under the assumptions (\ref{sum}), for each $u\in\R$ and $y\in S^\La$, there exists a unique cadlag function $Y^{u,y}\cn(-\infty,u]\to S^\La$ with final state $Y^{u,y}_u=y$ that satisfies the evolution equation
\be\label{duevol}
Y^{u,y}_{t-}:=\left\{\ba{ll}
\dis m(Y^{u,y}_t)\quad&\mbox{if }(m,t)\in\om,\\[5pt]
Y^{u,y}_t\quad&\mbox{otherwise,}\ea\right.\qquad(t\leq u).
\ee
Setting $\Yb_{u,t}(y):=Y^{u,y}_t$ $(u\geq t,\ y\in(S')^\La)$ then defines a collection $(\Yb_{u,t})_{u\geq t}$ of continuous maps from $(S')^\La$ into itself such that
\be\label{duflow}
\Yb_{s,s}=1\quand\Yb_{t,s}\circ\Yb_{u,t}=\Yb_{u,s}\qquad(u\geq t\geq s).
\ee
If $u\in\R$ and $Y_0$ is an $S^\La$-valued random variable with law $\mu$, independent of $\om$, then setting
\be\label{YY}
Y_t:=\Yb_{u,u-t}(Y_0)\qquad(t\geq 0)
\ee
defines a Feller process $(Y_t)_{t\geq 0}$ with generator $G$ as in (\ref{Gf}), just as in the case of the forward stochastic flow. The only difference is that, somewhat unusually, because of time-reversal, the sample paths of the Markov process defined in (\ref{YY}) are \emph{caglad}, that is, left-continuous with right limits.

\subsection{Additive duality}

In this subsection, we discuss additive interacting particle systems. We will see that each additive interacting particle system has a dual, which is another additive interacting particle system that we can think of as ``running backwards in time''.

We recall that a \emph{lattice} is a partially ordered set $S$ in which each two elements $a,b\in S$ have a greatest lower bound $a\wedge b$ and least upper bound $a\vee b$, called the \emph{infimum} and \emph{supremum} of $a$ and $b$. Each finite lattice $S$ contains a least element $0$ and greatest element $\top$. If $S$ is a finite lattice and $\La$ is countable, then the set $S^\La$, equipped with the product order, is a \emph{complete} lattice, meaning that each subset $A\sub S^\La$ has a greatest lower bound $\inf A$ and least upper bound $\sup A$. For a sequence $(x_n)_{n\in\N}$ in $S^\La$ we adopt the notation $\bigwedge_nx_n:=\inf\{x_n:n\in\N\}$ and $\bigvee_nx_n:=\sup\{x_n:n\in\N\}$. These are the pointwise infimum and supremum of the functions $x_n$ and $y_n$. The space $S^\La$ contains a least element $\un 0$ and greatest element $\un\top$, which are the functions that are constantly equal to $0$ and $\top$, respectively. By definition, a map $m\cn S^\La\to S^\La$ is \emph{additive} if 
\be
m(\un 0)=\un 0\quand m(x\vee y)=m(x)\vee m(y)\qquad(x,y\in S^\La),
\ee
and \emph{countably additive} if moreover $m\big(\bigvee_nx_n\big)=\bigvee_nm(x_n)$ for each sequence $x_n$ in $S^\La$. It is easy to see that each continuous additive map is also countably additive. 

A \emph{dual} of a partially ordered set $S$ is a partially ordered set $S'=(S',\leq)$ together with a bijection $S\ni a\mapsto a'\in S'$ such that
\be
a\leq b\quad\mbox{if and only if}\quad a'\geq b'.
\ee
All duals of $S$ are naturally isomorphic to each other, and setting $a'':=a$, we can naturally identify $S$ with the dual of $S'$. It is useful to see some examples. Each partially ordered set $(S,\leq)$ is dual to the same set equipped with the reversed order $(S,\geq)$ via the identity map $a\mapsto a':=a$. To see another example, let $\Pc(\De)$ be the set of all subsets of some set $\De$ and let $S\sub\Pc(\De)$ be closed under unions and intersections. Then $S$, equipped with the order of set inclusion, is a lattice. Letting $a^{\rm c}:=\De\beh a$ denote the complement of a set $a\sub\De$, one has that $S$ is dual to $S':=\{a^{\rm c}:a\in S\}$ via the complement map $a\mapsto a^{\rm c}$.

Let $S$ be a finite lattice, let $S'$ be a dual of $S$, and let $\La$ be a countable set. We equip $S^\La$ and $(S')^\La$ with the product order. Then $(S')^\La$ is dual to $S^\La$ via the map $x\mapsto x'$ that is defined pointwise by $x'(i):=x(i)'$ $(i\in\La)$. We define $\phi:S^\La\times(S')^\La\to\{0,1\}$ by
\be\label{phidef}
\phi(x,y):=1_{\txt\{x\leq y'\}}=1_{\txt\{x'\geq y\}}
\qquad\big(x\in S^\La,\ y\in (S')^\La\big).
\ee
If $m\cn S^\La\to S^\La$ and $\hat m\cn(S')^\La\to(S')^\La$ are maps, then we say that $\hat m$ is \emph{dual} to $m$ with respect to the duality function $\phi$ if
\be\label{dumaps}
\phi\big(m(x),y\big)=\phi\big(x,\hat m(y)\big)
\qquad\big(x\in S^\La,\ y\in (S')^\La\big).
\ee
The following lemma will be proved in Subsection~\ref{S:addproof}.

\bl[Dual maps]
Let\label{L:dumaps} $m\cn S^\La\to S^\La$ be countably additive. Then $m$ has a unique dual map $\hat m\cn(S')^\La\to(S')^\La$ with respect to the duality function $\phi$ from (\ref{phidef}). The dual map $\hat m$ is also countably additive. If $m$ is a local map, then so is $\hat m$.
\el

We say that an interacting particle system with generator of the form (\ref{Gf}) is \emph{additive} if the local state space $S$ is a finite lattice and the local maps $m\in\Gi$ are all additive. Assume from now on that this is the case and that the rates satisfy (\ref{sum}). Let $\om$ be a Poisson point set on $\Gi\times\R$ with intensity $\rho$ as in (\ref{rho}) and let $(\Xb_{s,t})_{s\leq t}$ be the stochastic flow defined in terms of the graphical representation $\om$. We will see in a moment that the maps $\Xb_{s,t}$ are countably additive, and give an expression for their duals.

For each $m\in\Gi$ let $\hat m\cn(S')^\La\to(S')^\La$ denote the dual map as in Lemma~\ref{L:dumaps}, which by that lemma is also local. We call
\be\label{hatom}
\hat\om:=\big\{(\hat m,t):(m,t)\in\om\big\}
\ee
the \emph{dual graphical representation}. Setting $\hat\Gi:=\{\hat m:m\in\Gi\}$, we have that $\hat\om$ is a Poisson point set on $\hat\Gi\times\R$ with intensity of the form (\ref{rho}), where the rates $(r_m)_{m\in\Gi}$ are replaced by the dual rates $(\hat r_m)_{m\in\hat\Gi}$ defined as $\hat r_{\hat m}:=r_m$ $(m\in\Gi)$. The following theorem will be proved in Subsection~\ref{S:addproof}.

\bt[Dual of additive particle system]
Assume\label{T:dual} that the rates $(r_m)_{m\in\Gi}$ satisfy (\ref{sum}). Then so do the dual rates $(\hat r_m)_{m\in\hat\Gi}$. Let $(\Xb_{s,t})_{s\leq t}$ be the stochastic flow defined in terms of $\om$ and let $(\Yb_{t,s})_{t\geq s}$ be the backward stochastic flow defined in terms of $\hat\om$. Then almost surely, for each $s\leq t$, the maps $\Xb_{s,t}$ and $\Yb_{t,s}$ are countably additive, and the dual map associated with $\Xb_{s,t}$ is given by $\hat\Xb_{s,t}=\Yb_{t,s}$.
\et

\noi
\textbf{Remark} We recall that conditions (\ref{sum}) (i) and (ii) imply that $(\Xb_{s,t})_{s\leq t}$ is well-defined. In the proof of Theorem~\ref{T:dual}, we will see that the dual rates $(\hat r_m)_{m\in\hat\Gi}$ satisfy condition (\ref{sum})~(ii) if and only if the forward rates $(r_m)_{m\in\Gi}$ satisfy condition (\ref{sum})~(iii), which is why we assume this too. Since $\un 0$ is a fixed point of $(\Xb_{s,t})_{s\leq t}$, by (\ref{finpert}), condition (\ref{sum})~(iii) moreover implies that $(\Xb_{s,t})_{s\leq t}$ maps the space $\big\{x\in S^\La:\sum_{i\in\La}1_{\{x(i)\neq 0\}}<\infty\big\}$ of finite configurations into itself.\med

Fix $t\geq 0$ and let $X_0$ and $Y_0$ be random variables with values in $S^\La$ and $(S')^\La$, respectively, independent of each other and of the graphical representation $\om$. Then by (\ref{XX}) and (\ref{YY}), setting
\be
X_s:=\Xb_{0,s}(X_0)\quand Y_s:=\Yb_{t,t-s}(Y_0)\quad(s\geq 0)
\ee
defines interacting particle systems $(X_s)_{s\geq 0}$ and $(Y_s)_{s\geq 0}$ with generators $G$ and $\hat G$ given by
\be
Gf(x):=\sum_{m\in\Gi}r_m\big\{f\big(m(x)\big)-f\big(x\big)\big\},\ 
\hat Gf(y):=\sum_{m\in\Gi}r_m\big\{f\big(\hat m(y)\big)-f\big(y\big)\big\}.
\ee
The duality $\hat\Xb_{s,t}=\Yb_{t,s}$ and the flow property give
\be
\phi(X_s,Y_{t-s})=\phi\big(\Xb_{0,s}(X_0),\Yb_{t,s}(Y_0)\big)=\phi\big(\Xb_{0,t}(X_0),Y_0\big)\qquad(0\leq s\leq t)
\ee
showing that the function $[0,t]\ni s\mapsto\phi(X_s,Y_{t-s})$ is constant. In particular, its values at $s=0$ and $s=t$ are the same, which after taking expectations yields the well-known duality relation
\be
\E\big[\phi(X_0,Y_t)\big]=\E\big[\phi(X_t,Y_0)\big]\qquad(t\geq 0),
\ee
where $X_0$ is independent of $Y_t$ and $X_t$ is independent of $Y_0$.

\subsection{Percolation representations}\label{S:percol}

In the present subsection we assume that $S=\{0,1\}$, equipped with the natural total order, $\La$ is a countable set, $\Gi$ is a collection of local maps $m\cn S^\La\to S^\La$, and $(r_m)_{m\in\Gi}$ are nonnegative rates satisfying (\ref{sum}). Let $\om$ be a Poisson point set on $\Gi\times\R$ with intensity as in (\ref{rho}), called the graphical representation, let $\hat\om$ be the dual graphical representation defined in (\ref{hatom}) and let $(\Xb_{s,t})_{s\leq t}$ and $(\Yb_{t,s})_{t\geq s}$ denote the stochastic flow and backward stochastic flow defined in terms of $\om$ and $\hat\om$. We will show that the stochastic flow $(\Xb_{s,t})_{s\leq t}$ and its dual $(\Yb_{t,s})_{t\geq s}$ have a percolation representation in terms of open paths. Our ultimate aim is to extend this percolation representation to the case that $S$ is a distributive lattice, but for the moment we recall the construction, which goes back to \cite{Har78,Gri79}, in the special case that $S=\{0,1\}$.

Let $\Pc(\La)$ denote the space of all subsets of $\La$, equipped with the order of set inclusion. Identifying a set with its indicator function, we have the natural isomorphism $\Pc(\La)\cong\{0,1\}^\La$. The partially ordered set $\Pc(\La)$ is dual to itself via the complement map $\Pc(\La)\ni x\mapsto x^{\rm c}:=\La\beh x$. The duality function $\phi$ from (\ref{phidef}) now takes the form
\be
\phi(x,y)=1_{\txt\{x\cap y=\emptyset\}}\qquad\big(x,y\in\Pc(\La)\big).
\ee
In the present setting, it is often more convenient to work with
\be\label{psidef}
\psi(x,y):=1-\phi(x,y)=1_{\txt\{x\cap y\neq\emptyset\}}\qquad\big(x,y\in\Pc(\La)\big).
\ee
Lemma~\ref{L:dumaps} then tells us that each countably additive map $m\cn\Pc(\La)\to\Pc(\La)$ has a unique dual map $\hat m$ such that
\be
\psi\big(m(x),y\big)=\psi\big(x,\hat m(y)\big)\qquad\big(x,y\in\Pc(\La)\big).
\ee
Countably additive maps on $\Pc(\La)$ have a simple structure. The following lemma will be proved in Subsection~\ref{S:percproof}.

\bl[Additive maps on sets]
For\label{L:mset} each $M\sub\La\times\La$, setting
\be\label{mM}
m(x):=\big\{j\in\La:\exists i\in x\mbox{ s.t.\ }(i,j)\in M\big\}\qquad\big(x\in\Pc(\La)\big)
\ee
defines a countably additive map $m\cn\Pc(\La)\to\Pc(\La)$. Conversely, each countably additive map $m\cn\Pc(\La)\to\Pc(\La)$ is of this form, with
\be\label{Mm}
M=\big\{(i,j)\in\La^2:j\in m(\{i\})\big\}.
\ee
Let $\hat m$ be the dual of $m$ with respect to the function $\psi$ from (\ref{psidef}) and let $\hat M$ be its corresponding set as in (\ref{Mm}). Then
\be\label{Madj}
\hat M=\big\{(j,i):(i,j)\in M\big\}.
\ee
\el

We call the set $M$ from (\ref{Mm}) the \emph{connection set} of the countably additive map $m$. We will give an explicit description of the connection sets of the maps $\Xb_{s,t}$ and $\Yb_{t,s}$. To this aim, we picture the graphical representation $\om$ by drawing space $\La$ on the horizontal axis, time $\R$ on the vertical axis, and then drawing for each $(m,t)\in\om$:
\begin{enumerate}
\item an arrow from $(i,t)$ to $(j,t)$ for each $i,j\in\La$ such that $i\neq j$ and $(i,j)\in M$, 
\item a blocking symbol \blok at $(i,t)$ for each $i\in\La$ such that $(i,i)\not\in M$.
\end{enumerate}
See Figure~\ref{fig:percol} for an example. For any $i,j\in\La$ and $s<t$, by definition, an \emph{open path} from $(i,s)$ to $(j,u)$ is a cadlag function $\ga\cn[s,u]\to\La$ such that $\ga_s=i$, $\ga_u=j$, and
\begin{enumerate}
\item if $\ga_{t-}\neq\ga_t$ for some $t\in(s,u]$, then there is an arrow from $(\ga_{t-},t)$ to $(\ga_t,t)$,
\item if $\ga_{t-}=\ga_t$ for some $t\in(s,u]$, then there is no blocking symbol at $(\ga_t,t)$.
\end{enumerate}
We write $(i,s)\leadsto(j,t)$ if there is an open path from $(i,s)$ to $(j,u)$. The following lemma will be proved in Subsection~\ref{S:percproof}.

\bl[Percolation representation]
For\label{L:perc} each $s\leq t$, let $M_{s,t}$ denote the connection set of the countably additive map $\Xb_{s,t}$ and let $\hat M_{t,s}$ denote the connection set of the countably additive map $\Yb_{t,s}$. Then
\be
M_{s,t}=\big\{(i,j)\in\La:(i,s)\leadsto(j,t)\big\}
\quand
\hat M_{t,s}=\big\{(j,i)\in\La:(i,s)\leadsto(j,t)\big\}.
\ee
\el

More explicitly, Lemma~\ref{L:perc} says that
\be\ba{r@{\ }r@{\,}c@{\,}l}\label{perc}
{\rm(i)}&\dis\Xb_{s,t}(x)&=&\dis\big\{j\in\La:\exists i\in x\mbox{ s.t.\ }(i,s)\leadsto(j,t)\big\},\\[5pt]
{\rm(ii)}&\dis\Yb_{t,s}(y)&=&\dis\big\{i\in\La:\exists j\in y\mbox{ s.t.\ }(i,s)\leadsto(j,t)\big\}.
\ec
This allows us to interpret the duality of $\Xb_{s,t}$ and $\Yb_{t,s}$ in terms of open paths, since
\be\label{pathdu}
\psi\big(\Xb_{s,t}(x),y\big)=1
\quad\desd\quad
\exists i\in x, j\in y\mbox{ s.t.\ }(i,s)\leadsto(j,t)
\quad\desd\quad
\psi\big(x,\Yb_{t,s}(y)\big)=1.
\ee

\subsection{Distributive lattices}

In the present subsection, we assume that $S$ is a finite distributive lattice, $\La$ is a countable set, $\Gi$ is a collection of local maps $m\cn S^\La\to S^\La$, and $(r_m)_{m\in\Gi}$ are nonnegative rates satisfying (\ref{sum}). We let $\om$ be a Poisson point set on $\Gi\times\R$ with intensity as in (\ref{rho}), define a dual graphical representation $\hat\om$ as in (\ref{hatom}), and let $(\Xb_{s,t})_{s\leq t}$ and $(\Yb_{t,s})_{t\geq s}$ denote the stochastic flow and backward stochastic flow defined in terms of $\om$ and $\hat\om$. Our aim is to show that also in this general setting, it is possible to give a percolation representation of $(\Xb_{s,t})_{s\leq t}$ and $(\Yb_{t,s})_{t\geq s}$ in terms of open paths.

We recall that a subset $A$ of a partially ordered set $\De$ is \emph{decreasing} (resp.\ \emph{increasing}) if $A\ni a\geq b$ implies $b\in A$ (resp.\ $A\ni a\leq b$ implies $b\in A$). We let $\Pc(\De)$ denote the set of all subsets of $\De$ and set
\bc
\dis\Pc^-(\De)&:=&\dis\big\{A\sub\De:A\mbox{ is decreasing}\big\},\\[5pt]
\dis\Pc^+(\De)&:=&\dis\big\{A\sub\De:A\mbox{ is increasing}\big\}.
\ec
By Birkhoff's representation theorem, a finite lattice $S$ is distributive if and only if there exists a finite partially ordered set $\De$ such that $S\cong\Pc^-(\De)$, where $\cong$ denotes isomorphism of partially ordered sets and we equip $\Pc^-(\De)$ with the order of set inclusion.

It follows that the state space of our interacting particle system satisfies
\be
S^\La\cong\Pc^-(\De)^\La\sub\Pc(\De)^\La\cong\{0,1\}^{\La\times\De}.
\ee
If we can extend the additive stochastic flow $(\Xb_{s,t})_{s\leq t}$ on $S^\La$ to an additive stochastic flow $(\ov\Xb_{s,t})_{s\leq t}$ on the larger space $\{0,1\}^{\La\times\De}$, then this will allow us to apply the theory from the previous subsection to construct a percolation representation for $(\Xb_{s,t})_{s\leq t}$ on the \emph{extended grid} $\La\times\De$. This basic idea can already be found in \cite{SS18}. In order to make it work, one needs to extend the additive local maps $m\cn S^\La\to S^\La$ that occur in the definition of the generator to additive local maps from the space $\{0,1\}^{\La\times\De}$ into itself. If the lattice $\La$ is finite, then it was shown in \cite[Lemma~13]{SS18} that this can always be done. Their construction does not work for infinite lattices, however. We therefore start by investivating how additive local maps that are initially defined on $S^\La$ can be extended to the larger space $\{0,1\}^{\La\times\De}$. We will especially be interested in extensions that have the property that their dual map also extends the dual of the original map.

We observe that $S=\Pc^-(\De)$ is dual to $S'\cong\Pc^+(\De)$ via the complement map $a\mapsto a^{\rm c}:=\De\beh a$. We equip the grid $\La$ with the \emph{trivial order}, in which $i\leq j$ if and only if $i=j$, and we equip the extended grid $\Lab:=\La\times\De$ with the product order. Elements of $\Lab$ are pairs $\ibf=(i,\sig)$ with $i\in\La$ and $\sig\in\De$. Then we have the natural isomorphisms of partially ordered sets
\be\label{hisms}
S^\La\cong\Pc^-(\Lab)\quand(S')^\La\cong\Pc^+(\Lab),
\ee
and these sets are dual via the complement map $\Pc^-(\Lab)\ni x\mapsto x^{\rm c}\in\Pc^+(\Lab)$. We view the local maps $m\in\Gi$ as maps on $\Pc^-(\Lab)$. To find a percolation representation of the stochastic flow $(\Xb_{s,t})_{s\leq t}$ and its dual $(\Yb_{t,s})_{t\geq s}$, it then suffices to extend each $m\in\Gi$ to a map $\mb\cn\Pc(\Lab)\to\Pc(\Lab)$ in such a way that if we set $\Gk:=\{\mb:m\in\Gi\}$ and $\rb_\mb:=r_m$, then the rates $(\rb_\mb)_{\mb\in\Gk}$ satisfy the conditions (\ref{sum}). Our main result is that this can always be done. We start by giving necessary and sufficient conditions for a countably additive map $\mb\cn\Pc(\Lab)\to\Pc(\Lab)$ to extend a countably additive map $m\cn\Pc^-(\Lab)\to\Pc^-(\Lab)$, and also for its dual $\hat\mb$ with respect to the duality function $\psi$ from (\ref{psidef}) to extend the dual $\hat m$ of $m$. The following lemma will be proved in Subsection~\ref{S:extproof}. Below, we write $\ibf<\jbf$ if $\ibf\leq\jbf$ and $\ibf\neq\jbf$. We write $\ibf,\ibf^\ast,\jbf,\jbf^\ast$ just to denote four arbitrary elements of $\Lab$, so below $\ibf^\ast$ is not some (hitherto undefined) function of $\ibf$.

\bl[Extensions of an additive map]
Let\label{L:addext} $m\cn\Pc^-(\Lab)\to\Pc^-(\Lab)$ be countably additive and let $N,N_\down$, and $N_\up$ be defined by
\bc\label{Nupdown}
\dis N&:=&\dis\big\{(\ibf,\jbf)\in\Lab^2:\jbf\in m(\{\ibf\}^\down)\big\}\quad\mbox{with}\quad\{\ibf\}^\down:=\big\{\ibf^\ast\in\Lab:\ibf^\ast\leq\ibf\big\},\\[5pt]
\dis N_\down&:=&\dis\big\{(\ibf,\jbf)\in N:\exists\ibf^\ast<\ibf\mbox{ s.t.\ }(\ibf^\ast,\jbf)\in N\big\},\\[5pt]
\dis N_\up&:=&\dis\big\{(\ibf,\jbf)\in N:\exists\jbf^\ast>\jbf\mbox{ s.t.\ }(\ibf,\jbf^\ast)\in N\big\}.
\ec
Then a countably additive map $\mb\cn\Pc(\Lab)\to\Pc(\Lab)$ satisfies $\mb(x)=m(x)$ for all $x\in\Pc^-(\Lab)$ if and only if its connection set $M$ satisfies $N\beh N_\down\sub M\sub N$. Moreover, one has $\hat\mb(y)=\hat m(y)$ for all $y\in\Pc^+(\Lab)$ if and only if $N\beh N_\up\sub M\sub N$.
\el

We equip the space of all countably additive maps $\mb\cn\Pc(\Lab)\to\Pc(\Lab)$ with a partial order such that $\mb_1\leq\mb_2$ if and only if $\mb_1(x)\sub\mb_2(x)$ for all $x\in\Pc(\Lab)$. It is easy to see that this is equivalent to the statement that the connection sets $M_1$ and $M_2$ of $\mb_1$ and $\mb_2$ satisfy $M_1\sub M_2$. We then see immediately from Lemma~\ref{L:addext} that the countably additive map $\mb\cn\Pc(\Lab)\to\Pc(\Lab)$ with connection set $M:=N$ is the maximal countably additive map that extends $m$. We call $\mb$ the \emph{maximal additive extension of $m$}. By Lemma~\ref{L:addext}, this map has the property that $\hat\mb$ extends $\hat m$. Lemma~\ref{L:addext} also tells us that the countably additive map $\mb\cn\Pc(\Lab)\to\Pc(\Lab)$ with connection set $M:=N\beh(N_\down\cap N_\up)$ is the minimal countably additive map such that $\mb$ extends $m$ and $\hat\mb$ extends $\hat m$. We call this map the \emph{minimal additive extension of $m$ and its dual}. Note that if we would set $M:=N\beh N_\down$, then $\mb$ would extend $m$ but $\hat\mb$ would not necessarily\footnote{A concrete example where this happens is the map $\mb^{2\to 1}_i$ defined in Subsection~\ref{S:twostage} below. As explained below (\ref{Nupdo}), this map can minimally be represented by a blocking symbol at $(i,1)$, but to also represent its dual one needs to add an arrow from $(i,1)$ to $(i,0)$.} extend $\hat m$, which is unwanted if we want to give a percolation representation not only for the stochastic flow $(\Xb_{s,t})_{s\leq t}$ but also for the dual stochastic flow $(\Yb_{t,s})_{t\geq s}$.

For finite grids $\La$, the maximal additive extension of a (countably) additive map $m$ has been introduced before in \cite[Lemma~13]{SS18}, but the definition of the minimal additive extension of $m$ and its dual seems to be new. Our final and main result is the following theorem that says that by extending all local maps using the minimal additive extensions of these maps and their duals, one can construct a percolation representation of any additive interacting particle system whose local state space is a distributive lattice and whose rates satisfy (\ref{sum}). A construction based on the maximal extensions of local maps would not work. Indeed, in most cases, these maximal extensions would fail to be local maps. The following theorem will be proved in Subsection~\ref{S:extproof}.

\bt[Minimal additive extension of a particle system and its dual]
Let\label{T:IPSext} $\Gi$ be a countable collection of additive local maps $m\cn\Pc^-(\Lab)\to\Pc^-(\Lab)$ and let $(r_m)_{m\in\Gi}$ be nonnegative rates satisfying (\ref{sum}). For each $m\in\Gi$, let $\mb$ denote the minimal additive extension of $m$ and its dual. Set $\Gk:=\{\mb:m\in\Gi\}$ and $\rb_\mb:=r_m$ $(m\in\Gi)$. Then each $\mb\in\Gk$ is a local map and the rates $(\rb_\mb)_{\mb\in\Gk}$ satisfy (\ref{sum}). Let $\omb:=\{(\mb,t):(m,t)\in\om\}$ and let $(\ov\Xb_{s,t})_{s\leq t}$ and $(\ov\Yb_{t,s})_{t\geq s}$ be defined in terms of $\omb$ and $\hat\omb:=\{(\hat\mb,t):(\mb,t)\in\omb\}$. Then almost surely, for each $s\leq t$, the map $\ov\Xb_{s,t}$ extends $\Xb_{s,t}$ and $\ov\Yb_{t,s}$ extends $\Yb_{t,s}$.
\et

\subsection{The two-stage contact process}\label{S:twostage}

Before we turn to the proofs, we demonstrate Theorem~\ref{T:IPSext} on a concrete example, the two-stage contact process introduced by Krone \cite{Kro99}. This is an interacting particle system with local state space $S=\{0,1,2\}$. To define the process and its dual, we need the local maps defined as
\be\ba{l}\label{mapi}
\dis m^{1\to 2}_i(x)(k):=\left\{\ba{@{}l@{\hspace{7pt}}l}
2&\mbox{if }k=i,\ x(i)=1,\\
x(k)&\mbox{otherwise,}\ea\right.\\[15pt]
\dis m^{1\to 0}_i(x)(k):=\left\{\ba{@{}l@{\hspace{5pt}}l}
0&\mbox{if }k=i,\ x(i)=1,\\
x(k)&\mbox{otherwise,}\ea\right.\\[15pt]
\dis m^{2\to 1}_i(x)(k):=\left\{\ba{@{}l@{\hspace{5pt}}l}
1&\mbox{if }k=i,\ x(i)=2,\\
x(k)&\mbox{otherwise,}\ea\right.\\[15pt]
\dis m^{\to 0}_i(x)(k):=\left\{\ba{@{}l@{\hspace{5pt}}l}
0&\mbox{if }k=i,\\
x(k)&\mbox{otherwise.}\ea\right.\\[15pt]
\dis m^\ast_{ij}(x)(k):=\left\{\ba{@{\,}l@{\hspace{10pt}}l}
1&\mbox{if }k=j,\ x(i)=2,\ x(j)=0\\
x(k)&\mbox{otherwise.}\ea\right.
\ec
The generator of the two-stage contact process then takes the form
\bc\label{Gdef}
\dis Gf(x)&:=&\dis\sum_{i,j\in\La}\la(i,j)\big\{f\big(m^\ast_{ij}(x)\big)-f\big(x\big)\big\}+\ga\sum_{i\in\La}\big\{f\big(m^{1\to 2}_i(x)\big)-f\big(x\big)\big\}\\[5pt]
&&\dis+\de_1\sum_{i\in\La}\big\{f\big(m^{1\to 0}_i(x)\big)-f\big(x\big)\big\}+\de_2\sum_{i\in\La}\big\{f\big(m^{\to 0}_i(x)\big)-f\big(x\big)\big\}.
\ec
Note that for the moment, we have not used the maps $m^{2\to 1}_i$. These will be needed only later, when we discuss the dual process. It is not hard to check that under the conditions
\be
\sup_{i\in\La}\sum_{j\in\La}\la(i,j)<\infty
\quand
\sup_{i\in\La}\sum_{j\in\La}\la(j,i)<\infty,
\ee
the rates of the process with generator $G$ satisfy (\ref{sum}), so the process is well-defined, can be constructed from a graphical representation $\om$, and the process started in a finite initial state stays finite. We can interpret the Markov process $(X_t)_{t\geq 0}$ with generator as in (\ref{Gdef}) as describing the time evolution of a forest. Here $X_t(i)=0,1,2$ means that the site $i$ is empty, occupied by a young tree, and occupied by an adult tree, respectively. An adult tree at $i$ gives birth with rate $\la(i,j)$ to a young tree at $j$, provided $j$ was previously empty. Because of the similarity to the contact process, we call $\la(i,j)$ the \emph{infection rate} from $i$ to $j$. Young trees mature at rate $\ga$, adult trees die at rate $\de_2$, and young trees die at rate $\de_1+\de_2$. The fact that young trees die at a higher rate than adult trees represents the biological reality that in many species, immature individuals have a higher death rate than adults.

The set $S=\{0,1,2\}$, equipped with the natural total order, is a distributive lattice. Setting $\De:=\{0,1\}$, we have that $S\cong\Pc^-(\De)$ via the map $0\mapsto\emptyset$, $1\mapsto\{0\}$, and $2\mapsto\{0,1\}$. Defining $\Lab:=\La\times\De$ and equipping it with an order as before, we have that
\be
S^\La\cong\Pc^-(\Lab)
\quand
(S')^\La\cong\Pc^+(\Lab).
\ee
It is straightforward to check that the local maps in (\ref{mapi}) are additive. If $m$ is one of the maps in (\ref{mapi}), then we let $\mb$ denote the minimal additive extension of $m$ and its dual. So $\mb^\ast_{ij}$ is the minimal additive extension of $m^\ast_{ij}$ and its dual, and so on. As in Theorem~\ref{T:IPSext}, we set $\omb:=\{(\mb,t):(m,t)\in\om\}$, we let $(\ov\Xb_{s,t})_{s\leq t}$ be the stochastic flow defined by $\omb$ and $(\ov\Yb_{t,s})_{t\geq s}$ the backward stochastic flow defined by $\hat\omb:=\{(\hat\mb,t):(\mb,t)\in\omb\}$. Then $(\ov\Xb_{s,t})_{s\leq t}$ and $(\ov\Yb_{t,s})_{t\geq s}$ are dual with respect to the duality function $\psi$ in (\ref{psidef}), and they extend the stochastic flows $(\Xb_{s,t})_{s\leq t}$ and $(\Yb_{t,s})_{t\geq s}$ of the two-type contact process and its dual, which are defined on the smaller spaces $\Pc^-(\Lab)$ and $\Pc^+(\Lab)$, respectively.

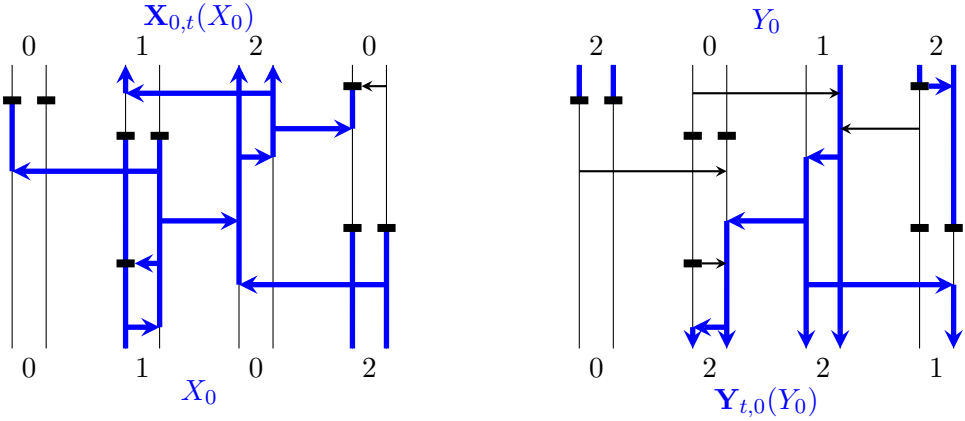
\begin{figure}[htb]
\begin{center}
\begin{tikzpicture}[xscale=1.5,yscale=1.25,>=stealth]
\begin{scope}[yscale=0.75]

\coordinate (a) at (2,0.3);
\coordinate (aa) at (2.3,0.3);
\coordinate (b) at (3,0.9);
\coordinate (c) at (4,0.9);
\coordinate (cc) at (4.3,0.9);
\coordinate (d) at (2,1.2);
\coordinate (dd) at (2.3,1.2);
\coordinate (e) at (2,1.8);
\coordinate (ee) at (2.3,1.8);
\coordinate (f) at (3,1.8);
\coordinate (g) at (2,2.5);
\coordinate (gg) at (2.3,2.5);
\coordinate (h) at (3,2.7);
\coordinate (hh) at (3.3,2.7);
\coordinate (i) at (2,3.0);
\coordinate (ii) at (2.3,3.0);
\coordinate (j) at (3,3.1);
\coordinate (jj) at (3.3,3.1);
\coordinate (k) at (4,3.1);
\coordinate (l) at (1,2.5);
\coordinate (m) at (1,3.5);
\coordinate (n) at (2,3.6);
\coordinate (o) at (3,3.6);
\coordinate (oo) at (3.3,3.6);
\coordinate (p) at (4,3.7);
\coordinate (q) at (4,1.7);
\coordinate (qq) at (4.3,1.7);

\draw[black!30,line width=3pt,rounded corners=5pt] (2,0)--(i);
\draw[black!30,line width=3pt,rounded corners=5pt] (2,0)--(a)--(aa)--(dd)--(d)--(i);
\draw[black!30,line width=3pt,rounded corners=5pt] (2,0)--(a)--(aa)--(ii);
\draw[black!30,line width=3pt,rounded corners=5pt] (2,0)--(a)--(aa)--(gg)--++(-1.3,0)--(m);
\draw[black!30,line width=3pt,rounded corners=5pt,->] (2,0)--(a)--(aa)--(ee)--++(0.7,0)--(3,4.1);
\draw[black!30,line width=3pt,rounded corners=5pt] (4,0)--(q);
\draw[black!30,line width=3pt,rounded corners=5pt] (4.3,0)--(qq);
\draw[black!30,line width=3pt,rounded corners=5pt] (4.3,0)--(cc)--++(-1.3,0)--(h)--(hh)--(jj)--++(0.7,0)--(p);
\draw[black!30,line width=3pt,rounded corners=5pt,->] (4.3,0)--(cc)--++(-1.3,0)--(h)--(hh)--(3.3,4.1);
\draw[black!30,line width=3pt,rounded corners=5pt,->] (4.3,0)--(cc)--++(-1.3,0)--(h)--(hh)--(oo)--++(-1.3,0)--(2,4.1);

\foreach \x in {1,...,4}
 \draw (\x,0) --++(0,4);
\foreach \x in {1,...,4}
 \draw (\x,0) ++(0.3,0) --++(0,4);

\draw[line width=3pt] (d)++(-0.08,0)--++(0.16,0);
\draw[line width=3pt] (i)++(-0.08,0)--++(0.16,0);
\draw[line width=3pt] (i)++(0.3,0)++(-0.08,0)--++(0.16,0);
\draw[line width=3pt] (m)++(-0.08,0)--++(0.16,0);
\draw[line width=3pt] (m)++(0.3,0)++(-0.08,0)--++(0.16,0);
\draw[line width=3pt] (p)++(-0.08,0)--++(0.16,0);
\draw[line width=3pt] (q)++(-0.08,0)--++(0.16,0);
\draw[line width=3pt] (q)++(0.3,0)++(-0.08,0)--++(0.16,0);

\draw[->,thick] (a)--++(0.3,0);
\draw[->,thick] (c)++(0.3,0)--(b);
\draw[->,thick] (d)++(0.3,0)--++(-0.22,0);
\draw[->,thick] (e)++(0.3,0)--++(0.7,0);
\draw[->,thick] (g)++(0.3,0)--(l);
\draw[->,thick] (h)--++(0.3,0);
\draw[->,thick] (j)++(0.3,0)--++(0.7,0);
\draw[->,thick] (o)++(0.3,0)--(n);
\draw[->,thick] (p)++(0.3,0)--++(-0.22,0);

\draw (1.15,0) node[below] {0};
\draw (2.15,0) node[below] {1};
\draw (3.15,0) node[below] {0};
\draw (4.15,0) node[below] {2};
\draw (1.15,4) node[above] {0};
\draw (2.15,4) node[above] {1};
\draw (3.15,4) node[above] {2};
\draw (4.15,4) node[above] {0};

\draw (2.65,-0.3) node[below] {$X_0$};
\draw (2.65,4.3) node[above] {$\Xb_{0,t}(X_0)$};
\end{scope}

\begin{scope}[xshift=5cm,yscale=0.75]

\coordinate (a) at (2,0.3);
\coordinate (aa) at (2.3,0.3);
\coordinate (b) at (3,0.9);
\coordinate (c) at (4,0.9);
\coordinate (cc) at (4.3,0.9);
\coordinate (d) at (2,1.2);
\coordinate (dd) at (2.3,1.2);
\coordinate (e) at (2,1.8);
\coordinate (ee) at (2.3,1.8);
\coordinate (f) at (3,1.8);
\coordinate (g) at (2,2.5);
\coordinate (gg) at (2.3,2.5);
\coordinate (h) at (3,2.7);
\coordinate (hh) at (3.3,2.7);
\coordinate (i) at (2,3.0);
\coordinate (ii) at (2.3,3.0);
\coordinate (j) at (3,3.1);
\coordinate (jj) at (3.3,3.1);
\coordinate (k) at (4,3.1);
\coordinate (l) at (1,2.5);
\coordinate (m) at (1,3.5);
\coordinate (mm) at (1.3,3.5);
\coordinate (n) at (2,3.6);
\coordinate (o) at (3,3.6);
\coordinate (oo) at (3.3,3.6);
\coordinate (p) at (4,3.7);
\coordinate (pp) at (4.3,3.7);
\coordinate (q) at (4,1.7);
\coordinate (qq) at (4.3,1.7);

\draw[black!30,line width=3pt,rounded corners=5pt] (1,4)--(m);
\draw[black!30,line width=3pt,rounded corners=5pt] (1.3,4)--(mm);
\draw[black!30,line width=3pt,rounded corners=5pt] (4,4)--(p)--(pp)--(qq);
\draw[black!30,line width=3pt,rounded corners=5pt] (4.3,4)--(qq);
\draw[black!30,line width=3pt,rounded corners=5pt,->] (3.3,4)--(3.3,-0.1);
\draw[black!30,line width=3pt,rounded corners=5pt,->] (3.3,4)--(hh)--(h)--(3,-0.1);
\draw[black!30,line width=3pt,rounded corners=5pt,->] (3.3,4)--(hh)--(h)--($(e)+(1,0)$)--(ee)--(2.3,-0.1);
\draw[black!30,line width=3pt,rounded corners=5pt,->] (3.3,4)--(hh)--(h)--($(e)+(1,0)$)--(ee)--(aa)--(a)--(2,-0.1);
\draw[black!30,line width=3pt,rounded corners=5pt,->] (3.3,4)--(hh)--(h)--($(c)+(-1,0)$)--++(1.3,0)--(4.3,-0.1);

\foreach \x in {1,...,4}
 \draw (\x,0) --++(0,4);
\foreach \x in {1,...,4}
 \draw (\x,0) ++(0.3,0) --++(0,4);

\draw[line width=3pt] (d)++(-0.08,0)--++(0.16,0);
\draw[line width=3pt] (i)++(-0.08,0)--++(0.16,0);
\draw[line width=3pt] (i)++(0.3,0)++(-0.08,0)--++(0.16,0);
\draw[line width=3pt] (m)++(-0.08,0)--++(0.16,0);
\draw[line width=3pt] (m)++(0.3,0)++(-0.08,0)--++(0.16,0);
\draw[line width=3pt] (p)++(-0.08,0)--++(0.16,0);
\draw[line width=3pt] (q)++(-0.08,0)--++(0.16,0);
\draw[line width=3pt] (q)++(0.3,0)++(-0.08,0)--++(0.16,0);

\draw[<-,thick] (a)--++(0.3,0);
\draw[<-,thick] (c)++(0.3,0)--(b);
\draw[<-,thick] (d)++(0.3,0)--++(-0.22,0);
\draw[<-,thick] (e)++(0.3,0)--++(0.7,0);
\draw[<-,thick] (g)++(0.3,0)--(l);
\draw[<-,thick] (h)--++(0.3,0);
\draw[<-,thick] (j)++(0.3,0)--++(0.7,0);
\draw[<-,thick] (o)++(0.3,0)--(n);
\draw[<-,thick] (p)++(0.3,0)--++(-0.22,0);

\draw (1.15,0) node[below] {0};
\draw (2.15,0) node[below] {2};
\draw (3.15,0) node[below] {2};
\draw (4.15,0) node[below] {1};
\draw (1.15,4) node[above] {2};
\draw (2.15,4) node[above] {0};
\draw (3.15,4) node[above] {1};
\draw (4.15,4) node[above] {2};

\draw (2.65,-0.4) node[below] {$\Yb_{t,0}(Y_0)$};
\draw (2.65,4.3) node[above] {$Y_0$};
\end{scope}
\end{tikzpicture}
\caption{Percolation representation of the two-stage contact process and its dual, the on-off contact process.}
\label{fig:percol}
\end{center}
\end{figure}

As in Subsection~\ref{S:percol}, we can represent each map $\mb$ occurring in the graphical representation $\omb$ in terms of arrows and blocking symbols, which then yields a representation of the stochastic flow $(\ov\Xb_{s,t})_{s\leq t}$ and the backward stochastic flow $(\ov\Yb_{t,s})_{t\geq s}$ in terms of open paths. We claim that by filling in the definitions (see Lemma~\ref{L:addext} and the discussion below it), one obtains that the maps we are interested in are represented in terms of arrows and blocking symbols as follows.
\[\ba{l@{\ }l}
\dis\mb^{1\to 2}_i&\mbox{an arrow from $(i,0)$ to $(i,1)$,}\\[5pt]
\dis\mb^{1\to 0}_i&\mbox{an arrow from $(i,1)$ to $(i,0)$ and a blocking symbol at $(i,0)$,}\\[5pt]
\dis\mb^{2\to 1}_i&\mbox{an arrow from $(i,0)$ to $(i,1)$ and a blocking symbol at $(i,0)$,}\\[5pt]
\dis\mb^{\to 0}_i&\mbox{blocking symbols at $(i,0)$ and $(i,1)$,}\\[5pt]
\dis\mb^\ast_{ij}&\mbox{an arrow from $(i,1)$ to $(j,0)$.}
\ea\]
See Figure~\ref{fig:percol} for an illustration. The map $\mb^{2\to 1}_i$ does not occur in the graphical representation for the process $X_t$, but in Figure~\ref{fig:percol} we see it occurring in the graphical representation for the dual process $Y_t$, if we moreover reverse the roles of the two elements of the set $\De$. This will be explained in more detail below.

It is instructive to see, step by step, how one arrives at the representations for the maps above. We only give the arguments in detail for the map $\mb^{2\to 1}_i$, which turns out to be the most interesting one. As already mentioned, we set $\De:=\{0,1\}$ so that $S=\{0,1,2\}$ is isomorphic to $\Pc^-(\De)$ via the map $0\mapsto\emptyset$, $1\mapsto\{0\}$, and $2\mapsto\{0,1\}$. Note that we can summarise this definition by saying that this is the map
\be
S\ni a\mapsto\big\{\sig\in\De:\sig<a\big\}\in\Pc^-(\De).
\ee
We equip $\La$ with the trivial order and $\Lab:=\La\times\De$ with the product order, which has the consequence that $S^\La\cong\Pc^-(\Lab)$ via the map
\be
S^\La\ni x\mapsto\big\{(\sig,i)\in\Lab:\sig<x(i)\big\}\in\Pc^-(\Lab).
\ee
In this identification, $m^{2\to 1}_i$ corresponds to the map from $\Pc^-(\Lab)$ into itself defined as
\be
m^{2\to 1}_i(x):=x\beh\{(i,1)\}\qquad\big(x\in\Pc^-(\Lab)\big).
\ee
To find the minimal additive extension of $m^{2\to 1}_i$ and its dual, we apply Lemma~\ref{L:addext}. For $j\in\La$ one has $\{(j,0)\}^\down=\{(j,0)\}$ and $\{(j,1)\}^\down=\{(j,0),(j,1)\}$. It follows that
\be
m^{2\to 1}_i\big(\{(i,1)\}^\down\big)=\{(i,0)\}
\quand
m^{2\to 1}_i\big(\{(j,\sig)\}^\down\big)=\{(j,\sig)\}^\down
\quad\mbox{in all other cases.}
\ee
As a consequence, for the map $m^{2\to 1}_i$, the set $N$ from (\ref{Nupdown}) is given by
\be
N=\big\{(\ibf,\jbf)\in\Lab^2:\jbf\leq\ibf\big\}\beh\big\{\big((i,1),(i,1)\big)\big\},
\ee
that is, $N$ consists of all pairs of the form $\big((j,\sig),(j,\tau)\big)$ with either $j\neq i$ and $\sig\geq\tau$, or $j=i$, $\sig$ is arbitrary, and $\tau=0$.
To obtain $N\beh N_\down$, we must remove from $N$ all pairs $\big((j,\sig),(j,\tau)\big)$ for which it is possible to make $\sig$ smaller while staying inside $N$. Similarly, $N\beh N_\up$ is obtained from $N$ by removing all elements for which it is possible to make $\tau$ larger while staying inside $N$. We see from this that
\bc\label{Nupdo}
\dis N\beh N_\down&=&\dis\big\{\big((j,\sig),(j,\sig)\big):j\neq i\big\}\cup\big\{\big((i,0),(i,0)\big)\big\},\\[5pt]
\dis N\beh N_\up&=&\dis\big\{\big((j,\sig),(j,\sig)\big):j\neq i\big\}\cup\big\{\big((i,0),(i,0)\big),\big((i,1),(i,0)\big)\big\}.
\ec
Following the discussion below Lemma~\ref{L:addext}, we see that the connection set $M$ of the minimal additive extension $\mb^{2\to 1}_i$ of $m^{2\to 1}_i$ and its dual is given by $M=N\beh(N_\down\cap N_\up)=(N\beh N_\down)\cup(N\beh N_\up)$, which in our case is equal to $N\beh N_\up$. To represent the map $\mb^{2\to 1}_i\cn\Pc(\Lab)\to\Pc(\Lab)$ with connection set $M=N\beh N_\up$ in terms of arrows and blocking symbols, we need to draw an arrow from $\ibf$ to $\jbf$ for each $(\ibf,\jbf)\in M$ with $\ibf\neq\jbf$, and to draw a blocking symbol at $\ibf$ for each $(\ibf,\ibf)\not\in M$. For the connection set $M=N\beh N_\up$, this amounts to drawing an arrow from $(i,1)$ to $(i,0)$, and drawing a blocking symbol at $(i,1)$. Note that if we would set $M:=N\beh N_\down$, then the arrow from $(i,1)$ to $(i,0)$ would be missing. The corresponding map $\mb\cn\Pc(\Lab)\to\Pc(\Lab)$ would still extend $m^{2\to 1}_i\cn\Pc^-(\Lab)\to\Pc^-(\Lab)$, but following the discussion below Lemma~\ref{L:addext}, we see that its dual $\hat\mb$ would not extend $\hat m^{2\to 1}_i$.

The map $\mb^{2\to 1}_i$ does not occur in the definition of the generator in (\ref{Gdef}) and is therefore also absent from the graphical representation $\omb$, but as already mentioned, we see something like it occurring in the dual graphical representation $\hat\omb$. By Lemma~\ref{L:mset}, if $\mb$ and $\hat\mb$ are countably additive maps from $\Pc(\Lab)$ into itself and $M$ and $\hat M$ are their connection sets, then $\hat M=\{(\jbf,\ibf):(\ibf,\jbf)\in M\}$. This means that we can find the percolation representation of $\hat\mb$ is given in terms of the percolation representation of $\mb$ by applying the principle:
\be
\mbox{reverse the arrows and keep the blocking symbols.}
\ee
Concretely, for our maps of interest, this means that we have the following representations (compare Figure~\ref{fig:percol}):
\[\ba{l@{\ }l}
\dis\hat\mb^{1\to 2}_i&\mbox{an arrow from $(i,1)$ to $(i,0)$,}\\[5pt]
\dis\hat\mb^{1\to 0}_i&\mbox{an arrow from $(i,0)$ to $(i,1)$ and a blocking symbol at $(i,0)$,}\\[5pt]
\dis\hat\mb^{2\to 1}_i&\mbox{an arrow from $(i,1)$ to $(i,0)$ and a blocking symbol at $(i,0)$,}\\[5pt]
\dis\hat\mb^{\to 0}_i&\mbox{blocking symbols at $(i,0)$ and $(i,1)$,}\\[5pt]
\dis\hat\mb^\ast_{ij}&\mbox{an arrow from $(j,0)$ to $(i,1)$.}
\ea\]
We observe that these maps are similar to the maps we had before, provided we exchange the roles of the two elements $0$ and $1$ of $\De$.

More formally, we observe that $\De$ is self-dual via the map $\sig\mapsto\sig':=1-\sig$. As a result, $\Lab=\La\times\De$ is self-dual via the map $(i,\sig)\mapsto(i,\sig')$, and the map
\be\label{xac}
\Pc^+(\Lab)\ni x\mapsto x'\in\Pc^-(\Lab)\quad\mbox{with}\quad x':=\{(i,\sig'):(i,\sig)\in x\}
\ee
is an isomorphism of partially ordered sets. It follows that $\Pc^-(\Lab)$ is self-dual via the map that is the concatenation of the complement map $\Pc^-(\Lab)\ni x\mapsto x^{\rm c}\in\Pc^+(\Lab)$ and the isomorphism $\Pc^+(\Lab)\ni x\mapsto x'\in\Pc^-(\Lab)$. More explicitly, this is the map
\be
\Pc^-(\Lab)\ni x\mapsto x^\ast\in\Pc^-(\Lab)\quad\mbox{with}\quad x^\ast:=\big\{(i,1-\sig):(i,\sig)\in\Lab\beh x\big\}.
\ee
More trivially, also $\Pc(\Lab)$ is self-dual via the map $x\mapsto x^\ast$. If we use this map, rather than the complement map, to define $\phi$ as in (\ref{phidef}) and then set $\psi:=1-\phi$ as before (compare (\ref{psidef})), then we obtain the function $\ti\psi\cn\Pc(\Lab)\times\Pc(\Lab)\to\{0,1\}$ given by
\be
\ti\psi(x,y)=1_{\txt\{x\cap y'\neq\emptyset\}}\qquad\big(x,y\in\Pc(\Lab)\big),
\ee
where $y'$ is defined in (\ref{xac}). Letting $\ti\mb$ denote the dual of $\mb$ with respect to $\ti\psi$, rather than the duality function $\psi$ in (\ref{psidef}), we see that
\be\ba{c}
\dis\ti\mb^{1\to 2}_i=\mb^{1\to 2}_i,\quad
\ti\mb^{1\to 0}_i=\mb^{2\to 1}_i,\quad
\ti\mb^{2\to 1}_i=\mb^{1\to 0}_i,\\[5pt]
\dis\ti\mb^{\to 0}_i=\mb^{\to 0}_i,\quad
\ti\mb^\ast_{ij}=\mb^\ast_{ji}.
\ec
Restricting to the smaller space $\Pc^-(\Lab)$ and translating back into the space $S^\La$, it follows that the two-stage contact process is dual to the interacting particle system with generator
\bc\label{tiGdef}
\dis\ti Gf(x)&:=&\dis\sum_{i,j\in\La}\la(j,i)\big\{f\big(m^\ast_{ij}(x)\big)-f\big(x\big)\big\}+\ga\sum_{i\in\La}\big\{f\big(m^{1\to 2}_i(x)\big)-f\big(x\big)\big\}\\[5pt]
&&\dis+\de_1\sum_{i\in\La}\big\{f\big(m^{2\to 1}_i(x)\big)-f\big(x\big)\big\}+\de_2\sum_{i\in\La}\big\{f\big(m^{\to 0}_i(x)\big)-f\big(x\big)\big\},
\ec
with respect to the duality function
\be
\ti\psi(x,y)=1_{\txt\{\exists i\in\La\mbox{ s.t.\ }x(i)+y(i)\geq 3\}}\qquad\big(x,y\in S^\La\big).
\ee
This duality was first observed by Krone \cite{Kro99}, who introduced the process. The generator in (\ref{tiGdef}) is almost the same as the one in (\ref{Gdef}), except that the infection rates $\la(i,j)$ have been replaced by the \emph{reversed infection rates} $\la^\dgg(i,j):=\la(i,j)$, and the map $m^{1\to 0}_i$ (which makes young trees die) has been replaced by $m^{2\to 1}_i$ (which turns adult trees into young trees). Krone interpreted this as individuals entering a non-reproductive state and called the process with the generator in (\ref{tiGdef}) the \emph{on-off contact process}.

The main contribution of the present paper is that we have shown that both the two-stage contact process and its dual, the on-off contact process, have a percolation representation that allows us to interpret the duality between these processes in terms of open paths as in (\ref{pathdu}). Moreover, we have shown that such a percolation representation of an additive interacting particle system and its dual can be found generally, and systematically, when the local state space is a distributive lattice.

For finite grids $\La$, this was already done in \cite{SS18}. In fact, Figure~2 of \cite{SS18} is almost identical to our present Figure~\ref{fig:percol}. The approach in \cite{SS18}, however, does not generalise to infinite grids since it is based on the maximal additive extensions of the local maps occurring in the generator. In general, the maximal additive extension of a local map need not be a local map, so that on infinite grids, the approach of \cite{SS18} will in many cases not yield a well-defined interacting particle system. By contrast, in Theorem~\ref{T:IPSext}, we have shown that by using the minimal additive extensions of local maps and their duals, one can always obtain a well-defined interacting particle system.

In this context, it is worth pointing out that Figure~2 of \cite{SS18} contains errors. It is suggested there that the maps shown in the figure are the maximal additive extensions of the local maps occurring in the generator of the two-stage contact process, but that is not the case. For example, the maximal additive extension of the map $m^{1\to 2}_i$ has the following percolation representation:
\be
\mbox{an arrow from $(i,0)$ to $(i,1)$ and for each $j\in\La$ an arrow from $(j,1)$ to $(j,0)$.}
\ee
Similarly, the percolation representations of the maximal additive extensions of the other maps occurring in (\ref{Gdef}) contain additional arrows from $(j,1)$ to $(j,0)$ for almost each grid point $j\in\La$. These arrows, that cause these maps to be nonlocal, are missing in Figure~2 of \cite{SS18}.

\section{Proofs}\label{S:proof}

\subsection{Additive duality}\label{S:addproof}

In this subsection we prove Lemma~\ref{L:dumaps} and Theorem~\ref{T:dual}. We need a bit of preparation. We recall that a function $f$ from one partially ordered space into another one is \emph{monotone} if $a\leq b$ implies $f(a)\leq f(b)$. It is easy to see that additive functions are monotone, since $a\leq b$ implies $f(a)\leq f(a)\vee f(b)=f(a\vee b)=f(b)$. In what follows, whenever $S$ is a finite lattice and $\La$ is a countable set, we equip the space $S^\La$ with the product order and product topology. A sequence $x_n\in S^\La$ is \emph{increasing} (resp.\ \emph{decreasing}) if $x_n\leq x_{n+1}$ (resp.\ $x_n\geq x_{n+1}$) for each $n$. We write $x_n\up x$ (resp.\ $x_n\down x$) if $x_n$ is increasing (resp.\ \emph{decreasing}) and $x_n\to x$. A function $f$ that is defined on $S^\La$ and takes values in a metrisable space is \emph{continuous with respect to increasing functions} if $x_n\up x$ implies $f(x_n)\to f(x)$. For each sequence $x_n\in S^\La$, we define
\bc
\dis\liminf_{n\to\infty}x_n&:=&\dis\sup\{y\in S^\La:\exists N\mbox{ s.t.\ }y\leq x_n\ \forall n\geq N\big\},\\[5pt]
\dis\limsup_{n\to\infty}x_n&:=&\dis\inf\{y\in S^\La:\exists N\mbox{ s.t.\ }y\geq x_n\ \forall n\geq N\big\}.
\ec
It is easy to see that $\liminf_{n\to\infty}x_n$ is the pointwise limit inferior of the functions $x_n$ and similarly $\limsup_{n\to\infty}x_n$ is the pointwise limit superior. As a consequence, $x_n\to x$ in the topology on $S^\La$ if and only if $\liminf_{n\to\infty}x_n=\limsup_{n\to\infty}x_n=x$. Let $S,T$ be finite lattices and let $\La,\De$ be countable. A function $f\cn S^\La\to T^\De$ is \emph{lower semi-continuous} if $x_n\to x$ implies $\liminf_{n\to\infty}f(x_n)\geq f(x)$. It is easy to see that a monotone function $f\cn S^\La\to S^\La$ is lower semi-continuous if and only if it is continuous with respect to increasing functions.
As a result, for an additive function $f\cn S^\La\to T^\De$, the following statements are equivalent:
\begin{enumerate}
\item $f$ is countably additive,
\item $f$ is continuous with respect to increasing functions,
\item $f$ is lower semi-continuous.
\end{enumerate}
We let $\Ci_{\rm add}(S^\La,T^\De)$ and $\Li_{\rm add}(S^\La,T^\De)$ denote the spaces of all additive functions $f\cn S^\La\to T^\De$ that are continous and lower semi-continuous, respectively. If $S$ is a finite lattice and $\La$ is a countable set, then we let
\be
S^\La_{\rm fin}:=\big\{x\in S^\La:\sum_{i\in\La}1_{\txt\{x(i)\neq 0\}}<\infty\big\}
\ee
denote the set of finite configurations in $S^\La$. For each $a\in S$ and $i\in\La$, we let $e^a_i\in S^\La$ be defined by $e^a_i(j):=a$ if $i=j$ and $:=0$ otherwise. Let $S'$ be a dual of $S$. We set $\psi:=1-\phi$, where $\phi$ is the duality function in (\ref{phidef}), that is,
\be\label{psidef2}
\psi(x,y):=1-\phi(x,y)\quad\mbox{with}\quad
\phi(x,y):=1_{\txt\{x\leq y'\}}
\qquad\big(x\in S^\La,\ y\in (S')^\La\big).
\ee
We will need the following lemma.

\bl[Additive indicator functions]
Let\label{L:addbin} $S$ be a finite lattice, let $S'$ be a dual of $S$, and let $\La$ be countable. Then
\be\ba{r@{\ }r@{\,}c@{\,}l}\label{addbin}
{\rm(i)}&\dis\Ci\big(S^\La,\{0,1\}\big)&=&\dis\big\{\psi(\,\cdot\,,y):y\in(S')^\La_{\rm fin}\big\},\\[5pt]
{\rm(ii)}&\dis\Li\big(S^\La,\{0,1\}\big)&=&\dis\big\{\psi(\,\cdot\,,y):y\in(S')^\La\big\}.
\ec
Moreover, $y\mapsto\psi(\,\cdot\,,y)$ is a bijection from $(S')^\La$ to $\Li\big(S^\La,\{0,1\}\big)$.
\el

\bpro
For each $y\in(S')^\La$, we have that $\psi(\un 0,y)=1_{\{\un 0\not\leq y'\}}=0$, while $\psi(\bigvee_nx_n,y)=1\ \desd\ \exists i\in\La, n\in\N\mbox{ s.t.\ }x_n(i)\not\leq y(i)\ \desd\ \exists n\in\N\mbox{ s.t.\ }\psi(x_n,y)=1$. From this, we see that $\psi(\,\cdot\,,y)$ is countably additive for each $y\in(S')^\La$, proving the inclusion $\supset$ in (\ref{addbin})~(ii).

Assume that $f\cn S^\La\to\{0,1\}$ is countably additive. For each $i\in\La$, define $A_i:=\{a\in S:f(e^a_i)=\un 0\}$. Since $f$ is additive we have $0\in A_i$, and $a,b\in A_i$ imply $a\vee b\in A_i$. It follows that $A_i$ has a unique greatest element $c_i$. Define $y\in(S')^\La$ by $y(i):=c'_i$. Then $f(e^a_i)=0\ \desd\ a\leq c_i\ \desd\ e^a_i\leq y'\ \desd\ \psi(e^a_i,y)=0$, which shows that $f(e^a_i)=\psi(e^a_i,y)$ for each $a\in S$ and $i\in\La$. Using moreover that $f$ and $\psi(\,\cdot\,,y)$ are countably additive, we conclude that $f=\psi(\,\cdot\,,y)$. This proves the inclusion $\sub$ in (\ref{addbin})~(ii).

Fix $i\in\La$. Then $y_1,y_2\in(S')^\La$ satisfy $\psi(e^a_i,y_1)=\psi(e^a_i,y_2)$ for all $a\in S$ if and only if $\{a\in S:a\leq y'_1(i)\}=\{a\in S:a\leq y'_2(i)\}$, which is equivalent to $y_1(i)=y_2(i)$. In particular, this implies that $y\mapsto\psi(\,\cdot\,,y)$ is a bijection from $(S')^\La$ to $\Li\big(S^\La,\{0,1\}\big)$.

If $f=\psi(\,\cdot\,,y)$ for some $y\in(S')^\La$, then $\Ri(f)$, as defined below (\ref{relev}), equals $\{i\in\La:y(i)\neq 0\}$. Moreover, if $x_1,x_2\in S^\La$ satisfy $x_1(i)=x_2(i)$ for all $i\in\Ri(f)$, then $f(x_1)=f(x_2)$. Using this and the characterisation of continuous maps $f\cn S^\La\to\{0,1\}$ given at the beginning of Subsection~\ref{S:IPS}, we see that $f$ is continuous if and only if $y\in(S')^\La_{\rm fin}$.
\epro

\bpro[Proof of Lemma~\ref{L:dumaps}]
For each $y\in(S')^\La$, the function $x\mapsto\psi\big(m(x),y\big)$ is countably additive, which by Lemma~\ref{L:addbin} implies that there exists a unique element $\hat m(y)\in(S')^\La$ such that $\psi\big(m(x),y\big)=\psi\big(x,\hat m(y)\big)$ for all $x\in S^\La$. It follows that $m$ has a unique dual map $\hat m$ with respect to the duality function $\psi$, or equivalently with respect to $\phi$. Since $\psi\big(x,\hat m(\un 0)\big)=\psi\big(m(x),\un 0\big)=0=\psi(x,\un 0)$ for all $x\in S^\La$, we see that $\hat m(\un 0)=\un 0$. Similarly, the fact that
\be
\psi\big(x,\hat m\big(\bigvee_ny_n\big)\big)=\psi\big(m(x),\bigvee_ny_n\big)=\bigvee_n\psi\big(m(x),y_n\big)=\bigvee_n\psi\big(x,\hat m(y_n)\big)=\psi\big(x,\bigvee_n\hat m(y_n)\big)
\ee
for each $x\in S^\La$ implies that $\hat m\big(\bigvee_ny_n\big)=\bigvee_n\hat m(y_n)$ for all $y_n\in(S')^\La$, proving that $\hat m$ is countably additive.

To complete the proof, we must show that if $m$ is local, then so is $\hat m$. It is easy to see (compare \cite[Exercise~2.33]{Swa26}) that a map $m\cn S^\La\to S^\La$ is local if and only if there exists a finite set $\Ga\sub\La$ and a map $n\cn S^\Ga\to S^\Ga$ such that
\be\label{local}
m(x)(i)=\left\{\ba{ll}
n\big((x(j))_{j\in\Ga}\big)(i)\quad&\mbox{if }i\in\Ga,\\[5pt]
x(i)\quad&\mbox{otherwise.}\ea\right.
\ee
Let $\Ga$ and $n$ be as in (\ref{local}), let $\hat n\cn(S')^\La\to(S')^\La$ be its dual and let $\hat m$ be defined in terms of $\hat n$ as in (\ref{local}). Then clearly $\hat m$ is local, and using the fact that
\be
\phi(x,y)=\prod_{i\in\La}1_{\txt\{x(i)\leq y'(i)\}},
\ee
it is straightforward to check that $\hat m$ is dual to $m$.
\epro

The proof of Theorem~\ref{T:dual} needs some further preparations. If $m\cn S^\La\to S^\La$ is countably additive, then for each $i,j\in\La$ we define $m[i,j]\cn S\to S$ by
\be\label{mat}
m[i,j](a):=m(e^a_i)(j)\qquad(a\in S,\ i,j\in\La).
\ee
Since $m$ is countably additive,
\be\label{mij}
m(x)(j)=\bigvee_{i\in\La}m[i,j]\big(x(i)\big).
\ee
We use similar notation for the dual map $\hat m$. It is then straightforward to check that
\be\label{mijdu}
\hat m[i,j]=\widehat{m[j,i]}\qquad(i,j\in\La),
\ee
where $\widehat{m[j,i]}$ is the dual of the map $m[i,j]\cn S\to S$. We let ${\rm id}\cn S\to S$ denote the identity map and we let $o\cn S\to S$ denote the zero map defined as $o(a):=0$ $(a\in S)$. We use similar notation for the identity map and zero map on $S'$. It is easy to see that $\hat{\rm id}={\rm id}$ and $\hat o=o$. We recall that for any map $m\cn S^\La\to S^\La$, the sets $\Di(m)$ and $\Ri(m)$ have been defined in (\ref{DR}).

\bl[Relevant lattice points]
Let\label{L:adrel} $m\cn S^\La\to S^\La$ be countably additive and let $\hat m\cn(S')^\La\to(S')^\La$ be its dual. Then
\be\label{adrel}
\Ri(m)=\big\{(i,j)\in\La^2:m(e^\top_i)(j)\neq 0\big\}=\big\{(i,j)\in\La^2:m[i,j]\neq o\big\},
\ee
and
\be\label{Dad}
\Di(m)=\big\{j\in\La:\exists i\neq j\mbox{ s.t.\ }m[i,j]\neq o\big\}
\cup\big\{j\in\La:m[j,j]\neq{\rm id}\big\}.
\ee
\el

\bpro
By definition $(i,j)\in\Ri(m)$ if there exist $x,y\in S^\La$ such that $x(k)=y(k)$ for all $k\neq i$ while $m(x)(j)\neq m(x)(j)$. Since $m$ is monotone, this is equivalent to the statement that there exist $x,y\in S^\La$ such that $x(i)=0$, $y(i)=\top$, and $x(k)=y(k)$ for all $k\neq i$, while $m(x)(j)\neq m(x)(j)$. We can then write
\be
x=\bigvee_{k\in\La\beh\{i\}}e^{x(k)}_k\quand y=x\vee e^\top_i
\ee
and use the fact that $m$ is countably additive to see that $(i,j)\in\Ri(m)$ if and only if $m(e^\top_i)(j)\neq 0$, proving the first equality in (\ref{adrel}). The second equality is an immediate consequence of this.

By definition, $\Di(m)$ is the set of all $j\in\La$ such that $m(x)(j)\neq x(j)$ for some $x\in S^\La$. If $m[i,j]=o$ for all $i\neq j$ and $m[j,j]={\rm id}$, then (\ref{mij}) implies that $j\not\in\Di(m)$. On the other hand, if $m[i,j]\neq o$ for some $i\neq j$, then setting $x:=e^\top_i$ we have $m(x)(j)\neq 0=x(j)$, while in the case that $m[j,j]\neq{\rm id}$ we can set $x:=e^a_j$ for some $a\in S$ to obtain $m(x)(j)\neq x(j)$ and hence $j\in\Di(m)$.
\epro

\bl[Summability condition]
Let\label{L:sumad} $S$ be a finite lattice, let $\La$ be countable, let $\Gi$ be a countable collection of additive local maps $m\cn S^\La\to S^\La$, and let $(r_m)_{m\in\Gi}$ be nonnegative rates. Then the summability condition (\ref{sum}) is equivalent to
\be\ba{l}\label{sumad}
\dis{\rm(i)}\ \sup_{i\in\La}\sum_{m\in\Gi}r_m1_{\txt\{m[i,i]\neq{\rm id}\}}<\infty,\\[5pt]
\dis{\rm(ii)}\ \sup_{i\in\La}\sum_{m\in\Gi}r_m\sum_{j\in\La\beh\{i\}}1_{\txt\{m[j,i]\neq o\}}<\infty,\\[5pt]
\dis{\rm(iii)}\ \sup_{i\in\La}\sum_{m\in\Gi}r_m\sum_{j\in\La\beh\{i\}}1_{\txt\{m[i,j]\neq o\}}<\infty.
\ec
\el

\bpro
By formula (\ref{adrel}) of Lemma~\ref{L:adrel}, conditions (ii) and (iii) of (\ref{sum}) are equivalent to conditions (ii) and (iii) of (\ref{sumad}), respectively. Let $\Di_1(m)$ and $\Di_2(m)$ denote the two sets on the right-hand side of formula (\ref{Dad}) of Lemma~\ref{L:adrel}. Then condition~(i) of (\ref{sumad}) is equivalent to
\be
\sup_{i\in\La}\sum_{m\in\Gi}r_m1_{\Di_2(m)}(i)<\infty,
\ee
which is certainly necessary for condition~(i) of (\ref{sum}) to hold. To see that it is also sufficient if we already know that condition~(ii) of (\ref{sum}), is satisfied, it suffices to note that
\be
1_{\Di_1(m)}(i)\leq\sum_{j\in\La\beh\{i\}}1_{\txt\{m[j,i]\neq o\}}.
\ee
\epro

\bpro[Proof of Theorem~\ref{T:dual}]
For the first statement of the theorem, by Lemma~\ref{L:sumad}, we may equivalently show that if the rates $(r_m)_{m\in\Gi}$ of the forward process satisfy (\ref{sumad}), then so do the dual rates $(\hat r_m)_{m\in\hat\Gi}$. By (\ref{mijdu}) and the fact that $\hat o=o$, we see that the dual rates $(\hat r_m)_{m\in\hat\Gi}$ satisfy (\ref{sumad})~(ii) if and only if the forward rates $(r_m)_{m\in\Gi}$ satisfy (\ref{sumad})~(iii). Likewise, (\ref{sumad})~(ii) for the forward rates is equivalent to (\ref{sumad})~(iii) for the dual rates, and using moreover that $\hat{\rm id}={\rm id}$ we see that the dual rates satisfy (\ref{sumad})~(i) if and only if the forward rates do.

Since the concatenation of two additive maps is additive, and since limits of additive maps are additive, by finite approximation, using \cite[Prop~4.23]{Swa26}, we find that $\Xb_{s,t}$ and $\Yb_{t,s}$ are additive for each $s\leq t$. By \cite[Thm~4.14]{Swa26} they are also continuous and hence countably additive. Since the duality function in (\ref{phidef}) is not continuous, we cannot prove duality by a straightforward approximation argument based on \cite[Prop~4.23]{Swa26} or the like. Instead, we argue as follows. We fix $y\in(S')^\La_{\rm fin}$ and $u\in\R$ and define for each $t\leq u$ a function $F_t\cn S^\La\to\{0,1\}$ by
\be\label{Fdef}
F_t(x):=\psi\big(\Xb_{t,u}(x),y\big)\qquad\big(t\leq u,\ x\in S^\La\big).
\ee
Since $\Xb_{t,u}$ is continuous and additive, we have $F_t\in\Ci_{\rm add}(S^\La,\{0,1\})$ for each $t\leq u$. By \cite[Thm~2.24 and Prop.~4.16]{Swa26}, the function $(-\infty,u]\ni t\mapsto F_t$ is piecewise constant and right-continuous, and solves the evolution equation
\be
F_{t-}=\left\{\ba{ll}
F_t\circ m\quad&\mbox{if }(m,t)\in\om,\\[5pt]
F_t\quad&\mbox{otherwise.}
\ea\right.
\ee
By Lemma~\ref{L:addbin}, for each $t\leq u$, there exists a unique $Y_t\in(S')^\La_{\rm fin}$ such that $F_t=\psi(\,\cdot\,,Y_t)$. It is then straightforward to check that $(-\infty,u]\ni t\mapsto Y_t$ solves the evolution equation
\be
Y_{t-}=\left\{\ba{ll}
\hat m(Y_t)\quad&\mbox{if }(m,t)\in\om,\\[5pt]
Y_t\quad&\mbox{otherwise.}
\ea\right.
\ee
This allows us to identity $Y_t=\Yb_{u,t}(y)$ $(t\leq u)$. Since this holds for general $y\in(S')^\La_{\rm fin}$, by (\ref{Fdef}) this proves that
\be
\psi\big(x,\Yb_{u,t}(y)\big)=\psi\big(\Xb_{t,u}(x),y\big)\qquad\big(t\leq u,\ x\in S^\La,\ y\in(S')^\La_{\rm fin}\big).
\ee
Using countable additivity, we can extend this to general $y\in(S')^\La$.
\epro

\subsection{Percolation representations}\label{S:percproof}

In this subsection, we prove Lemmas \ref{L:mset} and \ref{L:perc}.\med

\bpro[Proof of Lemma~\ref{L:mset}]
Identifying a set with its indicator function, we have the isomorphism of partially ordered sets $\Pc(\La)\cong\{0,1\}^\La$. In this picture, we can identify a countably additive map $m\cn\Pc(\La)\to\Pc(\La)$ with a countably additive map $m\cn\{0,1\}^\La\to\{0,1\}^\La$. Such a map is countably additive if and only if for each $j\in\La$, the coordinate map $m[j]\cn\{0,1\}^\La\to\{0,1\}$ defined as $m[j](x):=m(x)(j)$ is countably additive. We set $S:=\{0,1\}$, which is dual to itself via the map $a\mapsto a':=1-a$. Setting $S=\{0,1\}=S'$ in (\ref{psidef2}) yields the duality function
\be
\psi(x,y)=1_{\txt\{x\wedge y\neq\un 0\}}\qquad(x,y\in S^\La),
\ee
which corresponds to the function in (\ref{psidef}) if we identify $S^\La\cong\Pc(\La)$. By Lemma~\ref{L:addbin}, the countably additive maps $f\cn S^\La\to S$ are precisely the maps of the form $\psi(\,\cdot\,,y)$ with $y\in S^\La$. Thus $m\cn S^\La\to S^\La$ is countably additive if and only if for each $j\in\La$, there exists an $y_i\in S^\La$ such that $m[j](x)=\psi(x,y_j)$ for all $x\in S^\La$. Defining $Y\cn\La^2\to\{0,1\}$ by $Y(i,j):=y_j(i)$, this says that
\be\label{mY}
m(x)(j)=1\quad\desd\quad\exists i\in\La\mbox{ s.t.\ }x(i)=1\mbox{ and }Y(i,j)=1.
\ee
Translated into the set-valued language, this means that a map $m\cn\Pc(\La)\to\Pc(\La)$ is countably additive if and only if it is of the form (\ref{mM}), where the function $Y$ in (\ref{mY}) is the indicator function of the set $M$ in (\ref{mM}). Since $j\in m(\{i\})$ if and only if $(i,j)\in M$, we can recover $M$ from $m$ as in (\ref{Mm}). If $\hat m$ is the dual of $m$ with respect to the function $\psi$ in (\ref{psidef}), then $j\in m(\{i\})\ \desd\ \psi(m(\{i\}),\{j\})=1\ \desd\ i\in\hat m(\{j\})$, proving (\ref{Madj}).
\epro

\bpro[Proof of Lemma~\ref{L:perc}]
Defining $\Xb_{s,t}$ as in (\ref{perc})~(i), we see that for each $x\in\Pc(\La)$ and $s\in\R$, the function $[s,\infty)\ni t\mapsto X_t$ defined as $X_t:=\Xb_{s,t}(x)$ is cadlag and solves the evolution equation (\ref{evol}). This shows that if we take (\ref{perc})~(i) as an alternative definition of $\Xb_{s,t}$, then this definition is equivalent to our original definition; in other words (\ref{perc})~(i) holds. The proof of (\ref{perc})~(ii) is the same, and (\ref{perc}) implies Lemma~\ref{L:perc}.
\epro

\subsection{Extensions of additive maps}\label{S:extproof}

In this subsection, we prove Lemma~\ref{L:addext} and Theorem~\ref{T:IPSext}. Generalising notation introduced in (\ref{Nupdown}), for any $x\in\Pc(\Lab)$, we define $x^\down\in\Pc^-(\Lab)$ and $x^\up\in\Pc^+(\Lab)$ by
\be
x^\down:=\{\ibf^\ast\in\La:\exists\ibf\in x\mbox{ s.t.\ }\ibf^\ast\leq\ibf\}
\quand
x^\up:=\{\ibf^\ast\in\La:\exists\ibf\in x\mbox{ s.t.\ }\ibf\leq\ibf^\ast\}.
\ee
Note that $x$ is decreasing (resp.\ increasing) if and only if $x=x^\down$ (resp.\ $x=x^\up$).\med

\bpro[Proof of Lemma~\ref{L:addext}]
Let $m\cn\Pc^-(\Lab)\to\Pc^-(\Lab)$ and $\mb\cn\Pc(\Lab)\to\Pc(\Lab)$ be countably additive, let $N,N_\down$, and $N_\up$ be defined as in (\ref{Nupdown}), and let $M$ be the connection set of $\mb$. Our first aim is to show that $\mb(x)=m(x)$ for all $x\in\Pc^-(\Lab)$ if and only if
\be\label{ext}
N\beh N_\down\sub M\sub N.
\ee
We claim that (\ref{ext}) is equivalent to
\be\label{ext2}
M\sub N\quand\forall(\ibf,\jbf)\in N\ \exists\ibf^\ast\leq\ibf\mbox{ s.t.\ }(\ibf^\ast,\jbf)\in M.
\ee
Assume (\ref{ext2}). Then clearly $M\sub N$. Moreover, if $(\ibf,\jbf)\in N\beh N_\down$, then by (\ref{ext2}) we must have $(\ibf,\jbf)\in M$. This proves that (\ref{ext2}) implies (\ref{ext}). We next assume (\ref{ext}). Then clearly $M\sub N$. If $(\ibf,\jbf)\in N$, then either $(\ibf,\jbf)\in M$ or we can find $\ibf^\ast<\ibf$ such that $(\ibf^\ast,\jbf)\in N$. Then either $(\ibf^\ast,\jbf)\in M$ or we can find $\ibf^{\ast\ast}<\ibf$ such that $(\ibf^{\ast\ast},\jbf)\in N$. Continuing this process, using the definition of the partial order on $\Lab$ and the finiteness of $\De$, we find in finitely many steps a point $\ibf^\star\leq\ibf$ such that $(\ibf^\star,\jbf)\in M$. This completes the proof of the equivalence of (\ref{ext}) and (\ref{ext2}).

We next observe that
\be\label{Nmon}
(\ibf^\ast,\jbf)\in N\mbox{ and }\ibf^\ast\leq\ibf\quad\volgt\quad(\ibf,\jbf)\in N,
\ee
which follows directly from the definition of $N$ since $\{\ibf^\ast\}^\down\sub\{\ibf\}^\down$ and hence $m(\{\ibf^\ast\}^\down)\sub m(\{\ibf\}^\down)$ by the monotonicity of $m$. Using (\ref{Nmon}), we see that if $M$ satisfies (\ref{ext}) and hence (\ref{ext2}), then for each $\ibf,\jbf\in\Lab$,
\be
\jbf\in m(\{\ibf\}^\down)
\quad\desd\quad(\ibf,\jbf)\in N
\quad\desd\quad(\ibf^\ast,\jbf)\in M\mbox{ for some }\ibf^\ast\leq\ibf
\quad\desd\quad\jbf\in\mb(\{\ibf\}^\down).
\ee
This shows that $m(\{\ibf\}^\down)=\mb(\{\ibf\}^\down)$ for each $\ibf\in\Lab$. Since each element $x\in\Pc^-(\Lab)$ satisfies $x=x^\down$ and hence, by the definition of the latter, can be written as a countable union of sets of the form $\{\ibf\}^\down$, using the fact that $m$ and $\mb$ are countably additive, we see that $\mb(x)=m(x)$ for all $x\in\Pc^-(\Lab)$. This proves that (\ref{ext}) implies $\mb(x)=m(x)$ for all $x\in\Pc^-(\Lab)$.

Conversely, if $\mb(x)=m(x)$ for all $x\in\Pc^-(\Lab)$, then for each $\ibf,\jbf\in\La$,
\be
(\ibf,\jbf)\in M\ \volgt\ \jbf\in\mb(\{\ibf\})\ \volgt\ \jbf\in\mb\big(\{\ibf\}^\down\big)
\ \desd\ \jbf\in m\big(\{\ibf\}^\down\big)\ \desd\ (\ibf,\jbf)\in N,
\ee
which shows that $M\sub N$. Moreover, for each $\ibf,\jbf\in\La$,
\be
(\ibf,\jbf)\in N\ \volgt\ \jbf\in m\big(\{\ibf\}^\down\big)\ \desd\ \jbf\in\mb\big(\{\ibf\}^\down\big)\ \volgt\ \exists\ibf^\ast\leq \ibf\mbox{ s.t.\ }(\ibf^\ast,\jbf)\in M,
\ee
which shows that $M$ satisfies (\ref{ext2}) and hence (\ref{ext}). This completes the proof that $\mb(x)=m(x)$ for all $x\in\Pc^-(\Lab)$ if and only if (\ref{ext}) holds.

To complete the proof, we must show that $\hat\mb(y)=\hat m(y)$ for all $y\in\Pc^+(\Lab)$ if and only if $N\beh N_\up\sub M\sub N$. Letting $\hat M$ denote the connection set of $\hat\mb$ and applying what we have already proved to $\Lab$ equipped with the reversed order, we see that $\hat\mb(y)=\hat m(y)$ for all $y\in\Pc^+(\Lab)$ if and only if $\hat N\beh\hat N_\down\sub\hat M\sub\hat N$, where
\bc\label{NN}
\dis\hat N&:=&\dis\big\{(\jbf,\ibf)\in\Lab^2:\ibf\in\hat m(\{\jbf\}^\up)\big\},\\[5pt]
\dis\hat N_\down&:=&\dis\big\{(\jbf,\ibf)\in N:\exists\jbf<\jbf^\ast\mbox{ s.t.\ }(\jbf^\ast,\ibf)\in\hat N\big\}.
\ec
For any set $M\sub\Lab\times\Lab$, let us call $M^\dgg:=\{(\jbf,\ibf):(\ibf,\jbf)\in M\}$ the \emph{adjoint} of $M$. Then Lemma~\ref{L:mset} tells us that $\hat M=M^\dgg$. By definition, $\hat m$ maps increasing sets into increasing sets while $m$ maps decreasing sets into decreasing sets. Using this and duality, we see that
\be
\ibf\in\hat m(\{\jbf\}^\up)
\quad\desd\quad
\{\ibf\}^\down\cap\hat m(\{\jbf\}^\up)\neq\emptyset
\quad\desd\quad
m(\{\ibf\}^\down)\cap\{\jbf\}^\up\neq\emptyset
\quad\desd\quad
\jbf\in m(\{\ibf\}^\down),
\ee
which proves that $\hat N=N^\dgg$. It follows that
\be
(\hat N_\down)^\dgg=\big\{(\ibf,\jbf)\in N:\exists\jbf<\jbf^\ast\mbox{ s.t.\ }(\ibf,\jbf^\ast)\in N\big\}=N_\up.
\ee
Using this, by taking adjoints, we can rewrite our earlier condition $\hat N\beh\hat N_\down\sub\hat M\sub\hat N$ as $N\beh N_\up\sub M\sub N$.
\epro

The proof of Theorem~\ref{T:IPSext} needs a bit of preparation. In fact, the statement of the theorem needs a bit of clarification, as we now explain. Recall that $\Lab:=\La\times\De$ as defined above (\ref{hisms}). Setting $S:=\Pc^-(\De)$, $S':=\Pc^+(\De)$, and $T:=\Pc(\De)$, we have the natural identifications $\Pc^-(\Lab)\cong S^\La$, $\Pc^+(\Lab)\cong(S')^\La$, and $\Pc(\Lab)\cong T^\La$. In addition, we of course also have the isomorphism $\Pc(\Lab)\cong\{0,1\}^\Lab$. This has the consequence that the statement in Theorem~\ref{T:IPSext}, that says that the rates $(\rb_\mb)_{\mb\in\Gk}$ satisfy (\ref{sum}), can be interpreted in two ways: for each map $\mb\in\Gk$, we can interpret $\Di(\mb)$ and $\Ri(\mb)$ as subsets of $\La$ and $\La^2$, respectively, or alternatively as subsets of $\Lab$ and $\Lab^2$, depending on whether we view $\Pc(\De)$ as the local state space and $\La$ as the grid, or we view $\{0,1\}$ as the local state space and $\Lab=\La\times\De$ as the grid. We will show in both interpretations, the statement of Theorem~\ref{T:IPSext} is true.

In order to check (\ref{sum}), we will use its equivalent formulation (\ref{sumad}) in Lemma~\ref{L:sumad}. If we view $T$ as the local state space and $\La$ as the grid, then for each $\mb\in\Gk$ and $i,j\in\La$, we define $\mb[i,j]\cn T\to T$ as in (\ref{mat}). Similarly, if we view $\{0,1\}$ as the local state space and $\Lab$ as the grid, then for each $\ibf,\jbf\in\Lab$, we use (\ref{mat}) to define $\mb[\ibf,\jbf]\cn\{0,1\}\to\{0,1\}$. The only additive maps from $\{0,1\}$ into itself are the zero map and the identity map. It is straightforward to check that the connection set $M[i,j]$ of the additive map $\mb[i,j]\cn\Pc(\De)\to\Pc(\De)$ is given by
\be
M[i,j]=\big\{(\sig,\tau)\in\De^2:\mb\big[\big((i,\sig),(j,\tau)\big)\big]={\rm id}\big\}\qquad(i,j\in\La).
\ee
One has $\mb[i,j]={\rm id}$ if and only if $M[i,j]$ is the ``diagonal'' set $\{(\sig,\sig):\sig\in\De\}$, and $\mb[i,j]=o$ if and only if $M[i,j]=\emptyset$. As a consequence, if $\mb[(i,\sig),(i,\sig)]\neq{\rm id}$, then $\mb[i,i]\neq{\rm id}$, and if $\mb[(i,\sig),(j,\tau)]\neq o$ for some $(i,\sig)\neq(j,\tau)$, then either $i=j$ and $\mb[i,i]\neq{\rm id}$, or $i\neq j$ and $\mb[i,j]\neq o$. Using this, it is not hard to see that if the rates $(\rb_\mb)_{\mb\in\Gk}$ satisfy (\ref{sumad}) from the point of view that $T$ is the local state space and $\La$ is the grid, then they also satisfy (\ref{sumad}) in the other point of view, where $\{0,1\}$ is the local state space and $\Lab$ is the grid. In view of this, from now on, we view $T$ as the local state space and $\La$ as the grid, and set out to prove Theorem~\ref{T:IPSext} in this interpretation.

Recall that $S:=\Pc^-(\De)$ and $S':=\Pc^+(\De)$, and we have the natural identifications $\Pc^-(\Lab)\cong S^\La$ and $\Pc^+(\Lab)\cong(S')^\La$. We use this to define for each countably additive map $m\cn\Pc^-(\Lab)\to\Pc^-(\Lab)$ and $i,j\in\La$ a map $m[i,j]\cn S\to S$ as in (\ref{mat}), and similarly for countably additive maps on $\Pc^+(\Lab)$.

\bl[Nontrivial matrix elements]
Let\label{L:nontriv} $m\cn\Pc^-(\Lab)\to\Pc^-(\Lab)$ be countably additive and let $\mb\cn\Pc(\Lab)\to\Pc(\Lab)$ be the minimal additive extension of $m$ and its dual. Then for each $i,j\in\La$, one has $m[i,j]=o$ if and only if $\mb[i,j]=o$, and $m[i,j]={\rm id}$ if and only if $\mb[i,j]={\rm id}$. 
\el

\bpro
Let $M$ be the connection set of $\mb$ and let
\be
M[i,j]:=\big\{(\sig,\tau):\big((i,\sig),(j,\tau)\big)\in M\big\}\qquad(i,j\in\La).
\ee
Comparing the definition of a connection set in (\ref{Mm}) with the definition of $\mb[i,j]$ in (\ref{mat}), we see that $M[i,j]$ is the connection set of $\mb[i,j]$. Using this, Lemma~\ref{L:addext}, and the definition of the partial order on $\Lab$, we see that for each $i,j\in\La$, the map $\mb[i,j]$ is the minimal additive extension of $m[i,j]$ and its dual. Therefore, to complete the proof, it suffices to show that if $\ob$ and $\idb$ are the zero map and identity map on $T$, then I.\ $\ob$ is the minimal additive extension of $o$ and its dual, and II.\ $\idb$ is the minimal additive extension of ${\rm id}$ and its dual.

The proof of I.\ is easy: the duals of $o$ and $\ob$ are the zero maps on $S'$ and $T$, respectively, so $\ob$ extends $o$ and $\hat\ob$ extends $\hat o$. Since $\ob$ is the smallest map on $T$, this extension is obviously minimal. To also prove II., we use that if $m\cn\Pc^-(\De)\to\Pc^-(\De)$ is additive, and we set
\bc
\dis N&:=&\dis\big\{(\sig,\tau)\in\De^2:\tau\in m(\{\sig\}^\down)\big\},\\[5pt]
\dis N_\down&:=&\dis\big\{(\sig,\tau)\in N:\exists\sig^\ast<\sig\mbox{ s.t.\ }(\sig^\ast,\tau)\in N\big\},\\[5pt]
\dis N_\up&:=&\dis\dis\big\{(\sig,\tau)\in N:\exists\tau<\tau^\ast\mbox{ s.t.\ }(\sig,\tau^\ast)\in N\big\},
\ec
then the connection set $M$ of the minimal extension $\mb$ of $m$ and its dual is given by $M:=N\beh(N_\down\cap N_\up)$. Applying this to $m={\rm id}$, we see that
\be\ba{r}
\dis N=\big\{(\sig,\tau)\in\De^2:\tau\leq\sig\big\},\\[5pt]
\dis N_\down=N_\up=\big\{(\sig,\tau)\in\De^2:\tau<\sig\big\},
\ec
and $M=\{(\sig,\sig):\sig\in\De\}$. It follows that $\mb=\idb$, which is what we set out to prove.
\epro

\bpro[Proof of Theorem~\ref{T:IPSext}]
It is easy to see that a countably additive map $m\cn S^\La\to S^\La$ is local if and only if:
\begin{enumerate}
\item The set of all $(i,j)\in\La^2$ with $i\neq j$ for which $m[i,j]\neq o$ is finite.
\item The set of all $i\in\La$ for which $m[i,i]\neq{\rm id}$ is finite.
\end{enumerate}
An analogue statement holds for countably additive maps $\mb\cn T^\La\to T^\La$. In view of this and Lemma~\ref{L:nontriv}, if $m$ is an additive local map, then so is the minimal additive extension of $m$ and its dual. In the context of Theorem~\ref{T:IPSext}, this implies that the maps $\mb\in\Gk$ are all local. Moreover, combining Lemmas \ref{L:sumad} and \ref{L:nontriv}, we see that if the rates $(r_m)_{m\in\Gi}$ satisfy (\ref{sum}), then so do the rates $(\rb_\mb)_{\mb\in\Gk}$. It follows that the stochastic flow $(\ov\Xb_{s,t})_{s\leq t}$ and backward stochastic flow $(\ov\Yb_{t,s})_{t\geq s}$ are well-defined and given in terms of the unique solutions to evolution equations of the form (\ref{evol}) and (\ref{duevol}). Now if $x\in\Pc^-(\Lab)$ and $s\in\R$, then $t\mapsto X_{s,t}(x)$ solves the evolution equation that defines $(\ov\Xb_{s,t})_{s\leq t}$, which allows us to conclude that $\Xb_{s,t}(x)=\ov\Xb_{s,t}(x)$ for all $t\geq 0$. The proof that $\ov\Yb_{t,s}$ extends $\Yb_{t,s}$ is the same.
\epro


\begin{thebibliography}{Swa26}

\bibitem[Fox16]{Fox16}
E.~Foxall.
Duality and complete convergence for multi-type additive growth models.
\emph{Adv.\ Appl.\ Probab.}~48(1) (2016), 32--51.

\bibitem[Gri79]{Gri79}
D.~Griffeath.
\emph{Additive and Cancellative Interacting Particle Systems.}
Lecture Notes in Math.~724, Springer, Berlin, 1979.

\bibitem[Har78]{Har78}
T.E.~Harris.
Additive set-valued Markov processes and graphical methods.
\emph{Ann.\ Probab.}~6 (1978), 355--378.

\bibitem[Kro99]{Kro99}
S.~Krone.
The two-stage contact process.
\emph{Ann.\ Appl.\ Probab.}~9(2) (1999), 331--351.

\bibitem[SS18]{SS18}
A.~Sturm and J.M.~Swart.
Pathwise duals of monotone and additive Markov processes.
\emph{J.\ Theor.\ Probab.}~31(2) (2018), 932--983.

\bibitem[Swa26]{Swa26}
J.M.~Swart.
\emph{A Course in Interacting Particle Systems.}
Cambridge University Press, Cambridge, 2026.

\end{thebibliography}
\end{document}